\documentclass[11pt]{llncs}
\usepackage[letterpaper,hmargin=1.0in,vmargin=1.0150in]{geometry}
\usepackage{amsmath,amssymb,amsbsy,amsfonts,latexsym,amsopn,amstext,amsxtra,euscript,amscd,color}
\usepackage{array}
\usepackage{booktabs}
\usepackage{graphicx}
\usepackage{slashbox}
\usepackage{psfrag}
\usepackage{epsfig}

\def\00{{\bf 0}}
\def \+ {\oplus}

\def \\{\cr}
\def\({\left(}
\def\){\right)}

\newcommand{\BBR}{\mathbb{R}}

\pdfoutput=1

\providecommand{\newoperator}[3]{%
\newcommand*{#1}{\mathop{#2}#3}}
\newoperator{\FD}{\mathrm{FD}}{\nolimits}

\begin{document}
\title{Algorithms based on DQM with new sets of base functions for solving parabolic partial differential equations in $(2+1)$ dimension}
\author{Brajesh Kumar Singh\thanks{Address for Correspondence: Department of Applied Mathematics, School of Physical Sciences, Babasaheb Bhimrao Ambedkar University, Lucknow 226025 Uttar Pradesh INDIA}, Pramod Kumar}
\institute{Department of Applied Mathematics, School of Physical Sciences, \\Babasaheb Bhimrao Ambedkar University, Lucknow 226025 Uttar Pradesh INDIA\\
\email{bksingh0584@gmail.com}}
\date{\today}
\maketitle

\begin{abstract}
This paper deals with the numerical computations of two space dimensional time dependent \emph{parabolic partial differential equations} by adopting adopting an optimal five stage fourth-order strong stability preserving Runge - Kutta (SSP-RK54) scheme for time discretization, and three methods  of differential quadrature with different sets of modified B-splines as base functions, for space discretization: namely- $i)$ mECDQM: (DQM with modified extended cubic B-splines); $ii)$ mExp-DQM: DQM with modified exponential cubic B-splines, and $iii)$ MTB-DQM: DQM with modified trigonometric cubic B--splines. Specially, we implement these methods on  \emph{convection-diffusion} equation to convert them into a system of first order ordinary differential equations (ODEs), in time. The resulting system of ODEs can be solved using any time integration method, while we prefer SSP-RK54 scheme. All the three methods are found stable for two space convection-diffusion equation by employing matrix stability analysis method. The accuracy and validity of the methods are confirmed by three
test problems of two dimensional \emph{convection-diffusion} equation, which shows that the proposed approximate solutions by any of the method are in good agreement with the exact solutions.

\vspace{0.35cm}
\end{abstract}

\section{Introduction}
The convection-diffusion model can be expressed mathematically, which is a semi linear parabolic partial differential equation. Specially, we consider an initial value system of \emph{convection-diffusion} equation in $2$ dimension as:
\begin{equation}\label{eq-CDEs}
\left. \begin{split}
 & \frac{\partial u}{\partial  t} - \alpha_x  \frac{\partial^2 u}{\partial  x^2} -  \alpha_y  \frac{\partial^2 u}{\partial  y^2} + \beta_x \frac{\partial u}{\partial  x} + \beta_y \frac{\partial u}{\partial  y}=0, \\
&  u(x, y, 0) = u^0(x, y), \\
\end{split} \right\} in \quad \Omega \times (0, T],
\end{equation}
together with the Dirichlet boundary conditions: %$u(x, y, t) = f(x, y, t)$ ~~ on ~~ $\partial\Omega \times (0, T]$, i.e.,
\begin{equation}\label{Eq-BC}
\left.\begin{split}
& u(a, y, t) = f_{1}(y, t), \quad u(b, y, t) = f_{2}(y, t),\\
& u(x, c, t) = f_{3}(x, t), \quad u(x, d, t) = f_{4}(x, t),
\end{split}
\right\} on \quad  \partial\Omega \times (0, T],
\end{equation}
or Neumann boundary conditions:
\begin{equation}\label{Eq-NBC}
\left. \begin{split}
&\left.\frac{\partial u }{\partial x}\right|_{x=a, y} = g_{1}(y, t), \quad \left.\frac{\partial u }{\partial x}\right|_{x=b, y} = g_{2}(y, t),\\
&\left.\frac{\partial u }{\partial y}\right|_{x, y=c} = g_{3}(x, t), \quad \left.\frac{\partial u }{\partial y}\right|_{x, y=d} = g_{4}(x, t),
\end{split}
\right\}  on \quad \partial\Omega \times (0, T].
\end{equation}
where $\partial\Omega$ is the boundary of computational domain $\Omega=[a, b] \times [c, d]\subset \BBR^2$,
 $(0, T ]$ is time interval, and $f_{i}, g_{i}$ $(i=1,2,3,4)$ and $u_0$ are known smooth functions, and $u(x, y, t)$ denote heat or vorticity. The parameters: $\beta_x$ and $\beta_y$ are constant convective velocities while the constants $\alpha_x >0, \alpha_y>0$ are diffusion coefficients in the direction of $x$ and $y$, respectively.

The convection-diffusion models have remarkable applications in various branches of science and engineering, for instance, fluid motion, heat transfer, astrophysics, oceanography, meteorology, semiconductors, hydraulics, pollutant and sediment transport, and chemical engineering. Specially, in computational hydraulics and fluid dynamics to model convection-diffusion of quantities such as mass, heat, energy, vorticity \cite{R76}. Many researchers have paid their attention to develop some schemes which could produce accurate, stable and efficient solutions behavior of convection-diffusion problems, see  \cite{SPK16,AMS14,SS14,KA15,ED15,KA16} and the references therein.

In the last years, the convection-diffusion equation \eqref{eq-CDEs} has been solved numerically using various techniques: namely- finite element method \cite{DKK12}, Lattice Boltzmann method \cite{DCD93}, finite-difference scheme and higher-order compact finite difference schemes \cite{GM84,N91,S95,TG07,ZDC00}. A nine-point high-order compact implicit scheme proposed by Noye and Tan \cite{NT88} is third-order accurate in space and second-order accurate in time, and has a large zone of stability. An extension of higher order compact difference techniques for steady-state \cite{S95} to the time-dependent problems have been presented by Spotz and Carey \cite{SC01}, are fourth-order accurate in space and second or lower order accurate in time but conditionally stable. The fourth-order compact finite difference unconditionally stable scheme due to Dehghan and Mohebbi \cite{DM08} have the accuracy of order $O(h^4, t^4)$. A family of unconditionally stable finite difference schemes presented in \cite{CP05} have the accuracy of order $O(h^3, t^2)$. The schemes presented in \cite{KDD02} are based on high-order compact scheme and weighted time discretization, are second or lower order accurate in time and fourth-order accurate in space. The high-order alternating direction implicit (ADI) scheme with accuracy of order $O(h^4, t^2)$ proposed by Karaa and Zhang \cite{KZ04}, is unconditionally stable. A high-order unconditionally stable exponential scheme for unsteady $1$D convection-diffusion equation by Tian and Yua \cite{TY11} have the accuracy of order $O(h^4, t^2)$. A rational high-order compact alternating direction implicit (ADI) method have been developed for solving $2$D unsteady convection-diffusion problems \cite{T11} is unconditionally stable and have the accuracy of order $O(h^4, t^2)$. A unconditionally stable fourth-order compact finite difference approximation for discretizing spatial derivatives and the cubic $C_1$- spline collocation method in time, proposed by Mohebbi and Dehghan \cite{MD10}, have the accuracy of order $O(h^4, t^4)$.  An unconditionally stable, semi-discrete based on Pade approximation, by Ding and Zhang \cite{DZ09}, is fourth-order accurate in space and in time both. The most of schemes are based on the two-level finite difference approximations with Dirichlet conditions, and very few schemes have been developed to solve the convection-diffusion equation with Neumann's boundary conditions, see \cite{CLZF11,Y06} and references therein. The fourth-order compact finite difference scheme by Cao et al. \cite{CLZF11} is of $5$th-order accurate in time and 4th-order in the space. A high-order alternating direction implicit scheme based on fourth-order Pade approximation developed by You \cite{Y06} is unconditionally stable with the accuracy of order $O(h^4, t^2)$.

The differential quadrature method (DQM) dates back to Bellman et al. \cite{BKC72}. After the seminal paper of Bellman, various test functions have been proposed, among others, spline functions, sinc
function, Lagrange interpolation polynomials, radial base functions, modified cubic B-splines,
see \cite{K11,QC89,QC89a,SR92,S00,AS13,SA14,SAS16,SB16}, etc. Shu and Richards \cite{SR92} have generalized approach of DQM for numerical simulation of incompressible Navier-Stokes equation. The main goal of this paper is to find numerical solution of initial value system of $2$D \emph{convection-diffusion} equation with both kinds of boundary conditions (Dirichlet boundary conditions and Neumann boundary conditions), approximated by DQM with new sets of modified cubic B-splines (modified extended cubic B-splines, modified exponential cubic B-splines, modified trigonometric cubic B-splines) as base functions, and so called modified trigonometric cubic-B-spline differential quadrature method (MTB-DQM), modified exponential cubic-B-spline differential quadrature method (mExp-DQM) and third modified extended cubic-B-spline differential quadrature method (mECDQ). These methods are used to transform the convection diffusion problem into a system of first order ODEs, in time. The resulting system of ODEs can be solved by using various time integration algorithm, among them, we prefer SSP-RK54 scheme \cite{GKS09,SR02} due to its reduce storage space, which results in less accumulation errors. The accuracy and adaptability of the method is illustrated by three test problems of two dimensional convection diffusion equations.

The rest of the paper is organized into five more sections, which follow this introduction. Specifically, Section \ref{sec-metho-decr-tem} deals with the description of the methods: namely- MTB-DQM, mExp-DQM and mECDQ. Section  \ref{sec-impli} is devoted to the procedure for the implementation of describe above these methods for the system \eqref{eq-CDEs} together with the boundary conditions as in \eqref{Eq-BC} and \eqref{Eq-NBC}. Section \ref{sec-stab} deals with the stability analysis of the methods. Section \ref{sec-num} deals with the main goal of the paper is the numerical computation of three test problems. Finally, Section \ref{sec-conclu} concludes the results.

\section{Description of the methods} \label{sec-metho-decr-tem}
The differential quadrature method is an approximation to derivatives of a function is the weighted sum of the functional values at certain nodes \cite{BKC72}. The weighting coefficients of the derivatives is depend only on grids \cite{SR92}. This is the reason for taking the partitions $P[\Omega]$ of the problem domain $\Omega=\{ (x, y)\in R^2: 0\leq x, y\leq 1\}$ distributed uniformly as follows:
 $$P[\Omega] = \{ (x_i, y_j)\in \Omega: h_x=x_{i+1}-x_{i}, h_y=y_{j+1}-y_{j}, i\in \Delta_x, j\in \Delta_y\},$$
where $\Delta_x=\{1,2,\ldots,N_x\}, \Delta_y=\{1,2,\ldots,N_y\}$, and $
h_x=\frac{1}{N_x-1} \mbox{ and } h_y=\frac{1}{N_y-1}$
are the discretization steps in both $x$ and $y$ directions, respectively. That is, a uniform partition in each $x, y$-direction with the following grid points: $$0=x_1 < x_2< \ldots <x_i<\ldots <x_{N_x-1}< x_{N_x}=1,$$
 $$0 = y_1 < y_2< \ldots< y_j< \ldots< y_{N_y-  1}<y_{N_y}=1.$$
Let $(x_i, y_j)$ be the generic grid point and
$$u_{ij}\equiv u_{ij}(t)\equiv u(x_i, y_j, t),~~i\in \Delta_x, j\in \Delta_y.$$

The $r$-th order derivative of $u(x,y,t)$, for $r \in \{1, 2\}$, with respect to $x, y$ at $(x_i, y_j)$ for $i\in \Delta_{x}, j\in \Delta_{y}$ is approximated as follows:
\begin{equation}\label{eq-deri1-trig}
\begin{split}
& \left.\frac{\partial^r u}{\partial x^r}\right)_{ij} = \sum_{\ell=1}^{N_x} a_{i\ell}^{(r)} u_{\ell j}, \qquad  i\in\Delta_x,\quad i=1,2,\ldots N\\
& \left.\frac{\partial^r u}{\partial y^r}\right)_{ij} = \sum_{\ell=1}^{N_y} b_{j \ell}^{(r)} u_{i\ell}, \qquad  j\in\Delta_y,\quad i=1,2,\ldots M
\end{split}
\end{equation}
%\end{flushleft}
where the coefficients $a_{ip}^{(r)}$ and $b_{jp}^{(r)}$, the time dependent unknown quantities, are termed as the weighting functions of the $r$th-order derivative, to be computed using various type of base functions.

\subsection{The MTB-DQM} \label{sec-disc-trig}
\noindent The trigonometric cubic B-spline function $T_i=T_i(x)$ at node $i$ in $x$ direction, read as \cite{GV16,AMI14}:
\begin{eqnarray}\label{eq-cbs-trig}
T_i= \frac{1}{\omega} \left\{ \begin{array}{ll}
                 p^3(x_i),  & x \in [x_{i}, x_{i+1}) \\
                 p(x_i)\{p(x_i)q(x_{i+2})+p(x_{i+1})q(x_{i+3})\} + p^2(x_{i+1})q(x_{i+4}), &   x \in [x_{i+1}, x_{i+2})\\
                 q(x_{i+4})\{p(x_{i+1}q(x_{i+3})+ p(x_{i+2})q(x_{i+4})\} + p(x_{i})q^2(x_{i+3}),
                    &  x \in [x_{i+2}, x_{i + 3})\\
                  q^3(x_{i+4}),  &  x \in [x_{i + 3}, x_{i + 4})
       \end{array}  \right.
\end{eqnarray}
where $p(x_i)=\sin(\frac{x-x_i}{2});~~ q(x_i)=\sin(\frac{x_i-x}{2})$, $\omega = \sin\(\frac{h_x}{2}\) \sin\(h_x\) \sin\(\frac{3h_x}{2}\)$. Then the set $\{T_{0}, \ldots, T_{N_x}, T_{N_x+1}\}$ forms a base over the interval $[a, b]$. Setting
 \begin{equation*}
 \begin{split}
 & a_1 = \frac{\sin^2\(\frac{h_x}{2}\)}{ \sin(h_x) \sin\(\frac{3 h_x}{2}\)}; \quad a_2 = \frac{2}{ 1+ 2\cos(h_x)}; \quad a_4 = \frac{3}{ 4 \sin\(\frac{3 h_x}{2}\)}=-a_3;\\
 & a_5 = \frac{3 + 9 \cos(h_x)}{ 16 \sin^2\(\frac{h_x}{2}\) \(2 \cos\(\frac{h_x}{2}\) +  \cos\(\frac{3 h_x}{2}\)\) }; \quad a_6 =  \frac{ 3\cos^2(\frac{h_x}{2})}{ \sin^2\(\frac{h_x}{2}\) \(2 + 4 \cos\(h_x\)\) }.
 \end{split}
\end{equation*}
The values of $T_i$ and its first and second derivatives in the grid point $x_j$, denoted by $T_{ij}:=T_i(x_j)$, $T'_{ij}:=T'_i(x_j)$ and $T''_{ij}:=T''_i(x_j)$, respectively, read:
\begin{eqnarray}\label{tab-coeff-trig}
       T_{ij}= \left\{ \begin{array}{ll}
                 a_2,   &\mbox{ if }  i-j=0 \\
                 a_1,  & \mbox{ if } i-j=\pm 1\\
                 0,  & \mbox {otherwise}
       \end{array}\right.; \quad
              T'_{ij}= \left\{ \begin{array}{ll}
                  a_4,  & \mbox{ if }  i-j=1\\
                  a_3,  & \mbox{ if }  i-j=- 1\\
                 0, &\mbox {otherwise}
       \end{array}\right.; \mbox{\quad}
       T''_{ij}= \left\{ \begin{array}{ll}
                a_6,  &\mbox{ if }  i-j=0 \\
                a_5,  &\mbox{ if }  i-j=\pm 1\\
                 0  &  \mbox {otherwise}
       \end{array}\right.
\end{eqnarray}
The modified trigonometric cubic B-splines base functions are defined as follows \cite{AS13,SB16}:
\begin{flushleft}
\begin{equation}\label{eq-modi-CBs-trig}
\left\{ \begin{split}
& \tau_1 (x) = T_1(x) + 2 T_0(x)\\
&  \tau_2 (x) = T_2(x)  - T_0(x)\\
& \vdots\\
& \tau_j (x) = T_j(x), \mbox{ for } j = 3,4, \ldots, N_x-2\\
& \vdots\\
&  \tau_{N_x-1}(x) = T_{N_x-1}(x)  - T_{N_x+1}(x)\\
&  \tau_{N_x}(x) = T_{N_x}(x)  + 2 T_{N_x+1}(x)
\end{split}\right.
\end{equation}
\end{flushleft}
Now, the set $\{\tau_1, \tau_2,\ldots, \tau_{N_x}\}$ is the base over $[a, b]$. The procedure to define modified trigonometric cubic B-splines in $y$ direction, is followed analogously.

\subsubsection{The evaluation of the weighting coefficients $a_{ij}^{(r)}$ and $ b_{ij}^{(r)} (r=1,2)$}\label{sec-weight-trig}

\noindent In order to evaluate the weighting coefficients $a_{ip}^{(1)}$ of first order partial derivative in Eq. \eqref{eq-deri1-trig}, we use the modified trigonometric cubic B-spline $\tau_p(x)$, $p\in\Delta_x$ in DQ method as base functions. Setting $\tau'_{pi}:=\tau_p'(x_i)$ and $\tau_{p\ell}:=\tau_p(x_\ell)$. Using MTB-DQM, the approximate values of first order derivative is obtained as follows:
\begin{equation}\label{eq-deri1approxx-trig}
\tau'_{pi} = \sum_{\ell=1}^{N_x} a_{i\ell}^{(1)}  \tau_{p\ell}, \qquad  p,i \in \Delta_x.
\end{equation}
Setting  $\Im=[\tau_{p\ell}]$, $A=[a_{i\ell}^{(1)}]$, and $\Im'=[\tau'_{p\ell}]$. Eq. \eqref{eq-deri1approxx-trig} reduced to compact matrix form:
\begin{equation}\label{eqn-MTBS}
\Im A^T=\Im'.
%\Phi \overrightarrow{a}[i] = \overrightarrow{H}[i],  \mbox{ for } i= 1, \ldots, N_x,
\end{equation}
The coefficient matrix $\Im$ of order $N_x$ can be read from  \eqref{tab-coeff-trig} and \eqref{eq-modi-CBs-trig} as:
\begin{equation*}
 \Im= \left[
  \begin{array}{cccccccc}
  a_2+2a_1  $\quad$ & a_1 $\quad$&    $\quad$&           $\quad$&           $\quad$&         $\quad$&         \\
  $0    \quad$& a_2 $\quad$& a_1          $\quad$&           $\quad$&           $\quad$&         $\quad$&         \\
       $\quad$& a_1   $\quad$& a_2       $\quad$& a_1     $\quad$&           $\quad$&         $\quad$&         \\
       $\quad$&       $\quad$&  \ddots      $\quad$&   \ddots  $\quad$&   \ddots  $\quad$&         $\quad$&         \\
       $\quad$&       $\quad$&              $\quad$&   a_1     $\quad$&  a_2 $\quad$&    a_1  $\quad$&          \\
       $\quad$&       $\quad$&              $\quad$&           $\quad$&   a_1     $\quad$&   a_2$\quad$&    $0$   \\
       $\quad$&       $\quad$&              $\quad$&           $\quad$&           $\quad$&   a_1   $\quad$&   a_2+2a_1   \\
       \end{array}
\right]
\end{equation*}
and the columns of the matrix $\Im'$ read as:
\begin{equation*}
\Im'[1] = \left[ \begin{array}{c}
    2 a_4\\
    a_3-a_4\\
    0 \\
  \vdots  \\
          \\
  0     \\
  0     \\
\end{array}
\right],
\Im'[2] =\left[\begin{array}{c}
     a_4 \\
    $0$    \\
     a_3  \\
     $0$   \\
   \vdots  \\
           \\
    $0$    \\
\end{array}\right],
\ldots,
\Im'[N_x-1] =  \left[ \begin{array}{c}
     $0$     \\
     \vdots  \\
             \\
     $0$     \\
     a_4  \\
     $0$     \\
      a_3  \\
\end{array}
\right], \mbox{ and }
\Im'[N_x] =\left[ \begin{array}{c}
    $0$ \\
        \\
\vdots  \\
        \\
$0$     \\
  2a_4 \\
 a_3- a_4  \\
\end{array}
\right].
\end{equation*}

\subsection{The mExp-DQM }\label{sec-disc}
\noindent The exponential cubic B-splines function $\zeta_i=\zeta_i(x)$ at node $i$ in $x$ direction, reads \cite{KA15,ED15}:
\begin{eqnarray}\label{eq-cbs-exp}
\zeta_i= \frac{1}{h_x^3}\left\{ \begin{array}{ll}
                 b_2\{(x_{i-2}-x)-\frac{1}{p}\sinh(p(x_{i-2}-x))\},  & x \in [x_{i-2}, x_{i-1}) \\
                 a_1+b_1(x_{i}-x)+c_1\exp(p(x_{i}-x))+d_1\exp(p(x_{i}-x)), &   x \in [x_{i-1}, x_{i})\\
                 a_1+b_1(x-x_{i})+c_1\exp(p(x-x_{i}))+d_1\exp(p(x-x_{i})),
                    &  x \in [x_{i}, x_{i + 1})\\
                 b_2\{(x-x_{i+2})-\frac{1}{p}\sinh(p(x-x_{i+2}))\},  &  x \in [x_{i + 1}, x_{i + 2})\\
                 0,  &  \mbox {otherwise}
       \end{array}  \right.
\end{eqnarray}
where
 \begin{equation*}
 \begin{split}
 &a_1=\frac{pc h_x}{pch_x-s}; b_1=\frac{p}{2}\left( \frac{s^2-c(1-c)}{(pch_x-s)(1-c)}\right), b_2=\frac{p}{2(pch_x-s)}, c=\cosh(p h_x), s=\sinh(p h_x), \\
 &c_1=\frac{1}{4} \left\{\frac{\exp(-ph_x)(1-c)+s(\exp(-ph_x)-1) }{(pch_x-s)(1-c)}\right\}, d_1=\frac{1}{4} \left\{\frac{\exp(ph_x)(c-1)+s(\exp(ph_x)-1) }{(pch_x-s)(1-c)}
 \right\}.
 \end{split}
\end{equation*}
 The set $\{\zeta_{0}, \zeta_{1}, \zeta_{2}, \ldots, \zeta_{N_x}, \zeta_{N_x+1}\}$ forms a base over $[a, b]$. Setting the values of $\zeta_i$ and its first and second derivatives at $x_j$ by $\zeta_{ij}:=\zeta_i(x_j)$, $\zeta'_{ij}:=\zeta'_i(x_j)$ and $\zeta''_{ij}:=\zeta''_i(x_j)$, respectively. Then
 \begin{eqnarray}\label{tab-coeff-exp}
       \zeta_{ij}= \left\{ \begin{array}{ll}
                 1,   &\mbox{ if }  i-j=0 \\
               \frac{s-ph}{2(pc h_x-s)},  & \mbox{ if } i-j=\pm 1\\
                 0,  & \mbox {otherwise}
       \end{array}\right.; ~~
            \zeta_{ij}= \left\{ \begin{array}{ll}
                  - \frac{p(1-c)}{2(pc h_x - s)},  & \mbox{ if }  i-j=1\\
                   \frac{p(1-c)}{2(pc h_x - s)},  & \mbox{ if }  i-j=- 1\\
                 0, &\mbox {otherwise}
       \end{array}\right.; ~~
     \zeta_{ij}= \left\{ \begin{array}{ll}
                 -\frac{p^2 s}{(pc h_x - s)},  &\mbox{ if }  i-j=0 \\
                 \frac{p^2 s}{2(pc h_x - s)},  &\mbox{ if }  i-j=\pm 1\\
                 0  &  \mbox {otherwise}
       \end{array}\right.
\end{eqnarray}
The modified exponential cubic B-splines base functions are read as:
\begin{flushleft}
\begin{equation}\label{eq-modi-CBs-exp}
\left\{ \begin{split}
& \varsigma_1 (x) = \zeta_1(x) + 2\zeta_0(x)\\
&  \varsigma_2 (x) = \zeta_2(x)  - \zeta_0(x)\\
& \vdots\\
& \varsigma_j (x) = \zeta_j(x), \mbox{ for } j = 3,4, \ldots, N_x-2\\
& \vdots\\
&  \varsigma_{N_x-1}(x) = \zeta_{N_x-1}(x)  - \zeta_{N_x+1}(x)\\
&  \varsigma_{N_x}(x) = \zeta_{N_x}(x)  + 2 \zeta_{N_x+1}(x)
\end{split}\right.
\end{equation}
\end{flushleft}
The set $\{\varsigma_1, \varsigma_2,\ldots, \varsigma_{N_x}\}$ is a base over $[a, b]$. The procedure to define modified trigonometric cubic B-splines in $y$ direction, is followed analogously.

\subsubsection{{The evaluation of the weighting coefficients $a_{ij}^{(r)}$ and  $ b_{ij}^{(r)} (r=1,2)$}}
\label{subsubsec-weight-exp}

\noindent Setting $\varsigma'_{pi}:=\varsigma_p'(x_i)$ and $\varsigma_{p\ell}:=\varsigma_p(x_\ell)$ for all $p, \ell \Delta_x$. Using mExp-DQM, the approximate values of the first-order derivative is given by
\begin{equation}\label{eq-deri1approxx-exp}
\varsigma'_{pi} = \sum_{\ell=1}^{N_x} a_{i\ell}^{(1)}  \varsigma_{p\ell}, \qquad  p,i \in \Delta_x.
\end{equation}
Setting  $\beth=[\varsigma_{p\ell}]$, $A=[a_{i\ell}^{(1)}]$, and $\beth'=[\varsigma'_{p\ell}]$, then Eq. \eqref{eq-deri1approxx-exp} can be reduced to compact matrix form:
\begin{equation}\label{eq-mExp-S}
\beth A^T=\beth'.
%\Phi \overrightarrow{a}[i] = \overrightarrow{H}[i],  \mbox{ for } i= 1, \ldots, N_x,
\end{equation}
Let $\omega = \frac{p(1-c)h_x}{pc h_x-s}$ and $\theta =\frac{s-ph_x}{2(pc h_x-s)}$. Using Eqns. \eqref{tab-coeff-exp} and \eqref{eq-modi-CBs-exp}, the coefficient matrix $\beth$ of order $N_x$, read as:
\begin{equation*}
 \beth= \left[
  \begin{array}{cccccccc}
  \omega  $\quad$ & \theta $\quad$&    $\quad$&           $\quad$&           $\quad$&         $\quad$&         \\
  $0    \quad$& 1 $\quad$& \theta          $\quad$&           $\quad$&           $\quad$&         $\quad$&         \\
       $\quad$& \theta   $\quad$& 1       $\quad$& \theta      $\quad$&           $\quad$&         $\quad$&         \\
       $\quad$&       $\quad$&  \ddots      $\quad$&   \ddots  $\quad$&   \ddots  $\quad$&         $\quad$&         \\
       $\quad$&       $\quad$&              $\quad$&   \theta     $\quad$&   1 $\quad$&    \theta  $\quad$&          \\
       $\quad$&       $\quad$&              $\quad$&           $\quad$&   \theta     $\quad$&   1$\quad$&    $0$   \\
       $\quad$&       $\quad$&              $\quad$&           $\quad$&           $\quad$&   \theta   $\quad$&   \omega   \\
       \end{array}
\right]
\end{equation*}
and the columns of the matrix $\beth'$ read:
\begin{equation*}
\beth'[1] = \left[ \begin{array}{c}
    \omega/h_x\\
    -\omega/h_x\\
    0 \\
  \vdots  \\
          \\
  0     \\
  0     \\
\end{array}
\right],
\beth'[2] =\left[\begin{array}{c}
    \omega/2h_x \\
    $0$    \\
    -\omega/2h_x  \\
     $0$   \\
   \vdots  \\
           \\
    $0$    \\
\end{array}\right],
\ldots,
\beth'[N_x-1] =  \left[ \begin{array}{c}
     $0$     \\
     \vdots  \\
             \\
     $0$     \\
    \omega/2h_x  \\
     $0$     \\
     -\omega/2h_x   \\
\end{array}
\right], \mbox{ and }
\beth'[N_x] =\left[ \begin{array}{c}
    $0$ \\
        \\
\vdots  \\
        \\
$0$     \\
 \omega/h_x \\
-\omega/h_x   \\
\end{array}
\right].
\end{equation*}

\subsection{The mECDQ method} \label{sec-disc-ext}
\noindent The extended cubic B-splines function $\varphi_i=\varphi_i(x)$, in the $x$ direction and at the knots, reads \cite{KA16,SPK16,SPKa16}:
\begin{eqnarray}\label{eq-cbs-ext}
\varphi_i= \frac{1}{24}\left\{ \begin{array}{ll}
                 4(1-\lambda)P_{i-2}^{3}(x)+3 \lambda P_{i-2}^{4}(x),  & x \in [x_{i-2}, x_{i-1}) \\
                 24p+12 P_{i-1}^{}(x)+ 6(2+\lambda) P_{i-1}^{2}(x)-12 P_{i-1}^{3}(x)-3 \lambda P_{i-1}^{4}(x), &   x \in [x_{i-1}, x_{i})\\
                 24p - 12 P_{i+1}^{}(x)+ 6(2+\lambda) P_{i+1}^{2}(x)+12 P_{i+1}^{3}(x)-3 \lambda P_{i+1}^{4}(x),
                    &  x \in [x_{i}, x_{i + 1})\\
                 4(\lambda-1)P_{i+2}^{3}(x)+3 \lambda P_{i+2}^{4}(x),  &  x \in [x_{i + 1}, x_{i + 2})\\
                 0,  &  \mbox {otherwise}
       \end{array}  \right.
\end{eqnarray}
where $h_x P_i(x)=(x-x_i)$ and $24\wp = 4-\lambda$, and $\lambda$ is a free parameter \cite{SPK16}. The set $\{\varphi_{0}, \varphi_{1}, \varphi_{2}, \ldots, \varphi_{N_x}, \varphi_{N_x+1}\}$ forms a base over $[a, b]$. Let $12 h_x \theta=8+\lambda$ and $2 h_x^2 \omega =2+ \lambda$, then the values of $\varphi_i$ and its first and second derivatives in the grid point $x_j$, denoted by $\varphi_{ij}:=\varphi_i(x_j)$, $\varphi'_{ij}:=\varphi'_i(x_j)$ and $\varphi''_{ij}:=\varphi''_i(x_j)$, respectively,  read:

\begin{eqnarray}\label{tab-coeff-ext}
      \varphi_{ij}= \left\{ \begin{array}{ll}
                 \theta,   &\mbox{ if }  i-j=0 \\
                 \wp,  & \mbox{ if } i-j=\pm 1\\
                 0,  & \mbox {otherwise}
       \end{array}\right.; \quad
              2 h_x \varphi'_{ij}= \left\{ \begin{array}{ll}
                  1,  & \mbox{ if }  i-j=1\\
                  -1,  & \mbox{ if }  i-j=- 1\\
                 0, &\mbox {otherwise}
       \end{array}\right.; \mbox{\quad}
       \varphi''_{ij}= \left\{ \begin{array}{ll}
               -2 \omega,  &\mbox{ if }  i-j=0 \\
                \omega,   &\mbox{ if }  i-j=\pm 1\\
                 0  &  \mbox {otherwise}
       \end{array}\right.
\end{eqnarray}
The modified extended cubic B-splines base functions are defined as follows \cite{SPK16}:
\begin{flushleft}
\begin{equation}\label{eq-modi-CBs-ext}
\left\{ \begin{split}
& \phi_1 (x) = \varphi_1(x) + 2\varphi_0(x)\\
&  \phi_2 (x) = \varphi_2(x)  - \varphi_0(x)\\
& \vdots\\
& \phi_j (x) = \varphi_j(x), \mbox{ for } j = 3,4, \ldots, N_x-2\\
& \vdots\\
&  \phi_{N_x-1}(x) = \varphi_{N_x-1}(x)  - \varphi_{N_x+1}(x)\\
&  \phi_{N_x}(x) = \varphi_{N_x}(x)  + 2 \varphi_{N_x+1}(x)
\end{split}\right.
\end{equation}
\end{flushleft}
The set $\{\phi_1, \phi_2,\ldots, \phi_{N_x}\}$ is a base over $[a, b]$.

\subsubsection{\textbf{The evaluation of the weighting coefficients $a_{ij}^{(r)}$ and  $ b_{ij}^{(r)} (r=1,2)$}}

 Setting $\phi'_{pi}:=\phi_p'(x_i)$ and $\phi_{p\ell}:=\phi_p(x_\ell)$. Using mECDQ method, the approximate values of the first-order derivative is given by
\begin{equation}\label{eq-deri1approxx-ext}
\phi'_{pi} = \sum_{\ell=1}^{N_x} a_{i\ell}^{(1)}  \phi_{p\ell}, \qquad  p,i \in \Delta_x.
\end{equation}
Setting  $\Phi=[\phi_{p\ell}]$, $A=[a_{i\ell}^{(1)}]$, and $\Phi'=[\phi'_{pi}]$. Eq. \eqref{eq-deri1approxx-ext} can be re-written in compact matrix form as:
\begin{equation}\label{eqn-mEC-S}
\Phi A^T=\Phi'.
\end{equation}
Using Eqns. \eqref{tab-coeff-ext} and \eqref{eq-modi-CBs-ext},  the matrix $\Phi$ of order $N_x$ read as:
\begin{equation*}
 \Phi= \left[
  \begin{array}{cccccccc}
  $1$  $\quad$& \wp   $\quad$&              $\quad$&           $\quad$&           $\quad$&         $\quad$&         \\
  $0$  $\quad$& \theta$\quad$& \wp          $\quad$&           $\quad$&           $\quad$&         $\quad$&         \\
       $\quad$& \wp   $\quad$& \theta       $\quad$& \wp       $\quad$&           $\quad$&         $\quad$&         \\
       $\quad$&       $\quad$&  \ddots      $\quad$&   \ddots  $\quad$&   \ddots  $\quad$&         $\quad$&         \\
       $\quad$&       $\quad$&              $\quad$&   \wp     $\quad$&   \theta  $\quad$&    \wp  $\quad$&          \\
       $\quad$&       $\quad$&              $\quad$&           $\quad$&   \wp     $\quad$&   \theta$\quad$&    \wp   \\
       $\quad$&       $\quad$&              $\quad$&           $\quad$&           $\quad$&   \wp   $\quad$&   $1$    \\
       \end{array}
\right]
\end{equation*}
and the columns of the matrix $\Phi'$ read:
\begin{equation*}
\Phi'[1] = \left[ \begin{array}{c}
    -1/h_x\\
    1/h_x\\
    0 \\
  \vdots  \\
          \\
  0     \\
  0     \\
\end{array}
\right],
\Phi'[2] =\left[\begin{array}{c}
    -1/2h_x \\
    $0$    \\
    1/2h_x  \\
     $0$   \\
   \vdots  \\
           \\
    $0$    \\
\end{array}\right],
\ldots,
\Phi'[N_x-1] =  \left[ \begin{array}{c}
     $0$     \\
     \vdots  \\
             \\
     $0$     \\
    -1/2h_x  \\
     $0$     \\
     1/2h_x   \\
\end{array}
\right], \mbox{ and }
\Phi'[N_x] =\left[ \begin{array}{c}
    $0$ \\
        \\
\vdots  \\
        \\
$0$     \\
 -1/h_x \\
1/h_x   \\
\end{array}
\right].
\end{equation*}

\noindent Using ``Thomas algorithm" the system \eqref{eqn-MTBS}, \eqref{eq-mExp-S} and \eqref{eqn-mEC-S} have been solved for the weighting coefficients $a_{i1}^{(1)}, a_{i2}^{(1)}, \ldots, a_{i N_x}^{(1)}$, for all  $i \in \Delta_x$.

Similarly, the weighting coefficients $b_{i j}^{(1)}$, in either case, can be computed by employing these modified cubic B-splines in the $y$ direction.

Using $a_{i j}^{(1)}$ and $b_{i j}^{(1)}$, the weighting coefficients, $a_{i j}^{(r)}$ and $b_{i j}^{(r)}$ (for $r \geq 2$) can be computed using the Shu's recursive formulae \cite{S00}:
\begin{flushleft}
\begin{equation}\label{eq-coeff2}
\left\{ \begin{split}
& a_{i j}^{(r)} = r \left( a_{i j}^{(1)} a_{i i}^{(r-1)} - \frac{a_{i j}^{(r-1)}}{ x_i - x_j}\right), i\ne j: i, j\in \Delta_x,
\\&  a_{i i}^{(r)} = - \sum_{i = 1, i \ne j}^{N_x} a_{i j}^{(r)}, i= j: i, j\in \Delta_x.\\
& b_{i j}^{(r)} = r \left(b_{i j}^{(1)} b_{i i}^{(r-1)} - \frac{b_{i j}^{(r-1)}}{ y_i - y_j}\right), i\ne j: i, j \in \Delta_y
 \\& b_{i i}^{(r)} = - \sum_{i = 1, i \ne j}^{N_y} b_{i j}^{(r)}, i = j: i, j \in \Delta_y.\\
\end{split}\right.
\end{equation}
\end{flushleft}
where $r$ denote $r$-th order spatial derivative. In particular, the weighting coefficients $a_{i j}^{(2)}, b_{i j}^{(2)}$ of order $2$ can be obtained by taking $r=2$ in \eqref{eq-coeff2}.

\section{Implementation of the method for $2$D convection-diffusion equation}\label{sec-impli}
After computing the approximate values of first and second order spatial partial derivatives from one of the above three methods, one can re-write Eq \eqref{eq-CDEs} as follows:

 \begin{equation}\label{eq-CBEs-ode}
 \left\{ \begin{split}
& \frac{du_{ij}}{dt} = \alpha_x \sum_{k=1}^{N_x} a_{ik}^{(2)} u_{kj} + \alpha_y \sum_{k=1}^{N_y} b_{jk}^{(2)} u_{ik} - \beta_x \sum_{k=1}^{N_x} a_{ik}^{(1)} u_{kj} - \beta_y \sum_{k=1}^{N_y} b_{jk}^{(1)} u_{ik},\\
& u_{ij}(t=0)=u^0_{ij}.
\end{split}\right.
\end{equation}
In case of the Dirichlet conditions, the solutions on boundaries can directly read from \eqref{Eq-BC} as:
 \begin{equation}\label{eq-DBC-dis}
u_{1j} = f_{1}(y_j, t); ~~ u_{N_xj} = f_{2}(y_j, t); ~~ u_{i1} = f_{3}(x_i, t); u_{iN_y} = f_{4}(x_i, t), \forall ~~ i \in \Delta_x, j \in \Delta_y, t \in (0, T].
\end{equation}

\noindent On the other hand, if the boundary conditions are Neumann or mixed type, then the solutions at the boundary are obtained by using any above methods (MTB-DQM, mExp-DQM or mECDQ method) on the boundary, which gives a system of two equations. On solving it we get the desired solution on the boundary as follows:

From Eq. \eqref{eq-deri1-trig} with $r=1$ and the Neumann boundary conditions \eqref{Eq-NBC} at $x= a$ and $x=b$, we get
\begin{flushleft}
 \begin{equation}\label{eq-N1-dis}
\left. \begin{split}
& \sum_{k=1}^{N_x} a_{1k}^{(1)} u_{kj} = g_1(y_j,t),\\
& \sum_{k=1}^{N_x} a_{Nk}^{(1)} u_{kj} = g_2(y_j,t),
\end{split}\right\} \quad j \in \Delta_y.
\end{equation}
\end{flushleft}
In terms of matrix system for $u_{1j}, u_{N_x j}$, the above equation can be rewritten as
\begin{equation}\label{eq-N2-dis}
\left[\begin{array}{cc}
a_{11}^{(1)}  & a_{1N_x}^{(1)}  \\
a_{N_x 1}^{(1)} & a_{N_x N_y}^{(1)} \\
\end{array} \right] \left[
\begin{array}{c}
u_{1j} \\
u_{N_x j} \\
\end{array} \right]
=
\left[
\begin{array}{c}
S_j^a \\
S_j^b \\
\end{array} \right],
\end{equation}
where $S_j^a= g_1(y_j,t)-\sum_{k=2}^{N_x-1} a_{1k}^{(1)} u_{kj}$ and $ S_j^b= g_2(y_j,t)- \sum_{k=2}^{N_x-1} a_{N_x k}^{(1)} u_{kj}$. On solving \eqref{eq-N2-dis}, for the boundary values $u_{1j} $ and $ u_{N_xj}, j \in \Delta_y$, we get
\begin{equation}\label{eq-N3-dis}
u_{1j} = \frac{S_j^a a_{N_xN_x}^{(1)} - S_j^b a_{1N_x}^{(1)} }{a_{11}^{(1)}a_{N_xN_x}^{(1)}-a_{N_x1}^{(1)}a_{1N_x}^{(1)} },  \qquad
u_{N_x j} = \frac{S_j^b a_{11}^{(1)} - S_j^a a_{N_x1}^{(1)}}{a_{11}^{(1)}a_{N_xN_x}^{(1)}-a_{N_x1}^{(1)}a_{1N_x}^{(1)} }.
\end{equation}

Analogously, for the Neumann boundary conditions \eqref{Eq-NBC} at $y= c$ and $y=d$, the solutions for the boundary values $u_{i1} $ and $ u_{iN_y}, i \in \Delta_x$ can be obtained as:
\begin{equation}\label{eq-N2y-dis}
%\left\{
 \begin{split}
&u_{i1} = \frac{S_i^c b_{N_yN_y}^{(1)} - S_i^d b_{1N_y}^{(1)} }{b_{11}^{(1)}b_{N_y N_y}^{(1)}-b_{N_y 1}^{(1)}b_{1 N_y}^{(1)} }, \qquad  u_{i N_y} = \frac{S_i^d b_{11}^{(1)} - S_i^c b_{N_y 1}^{(1)}}{b_{11}^{(1)}b_{N_y N_y}^{(1)}-b_{N_y 1}^{(1)}b_{1 N_y}^{(1)} },
\end{split}
\end{equation}
where $S_i^c= g_3(x_i,t) - \sum_{k=2}^{N_y-1} b_{1k}^{(1)} u_{ik}$ and $ S_i^d= g_4(x_i,t) - \sum_{k=2}^{N_y-1} b_{N_y k}^{(1)} u_{ik}$.

After implementing the boundary values, Eq \eqref{eq-CBEs-ode} can be written in compact matrix form as follows:
\begin{equation}\label{eq-ode-IVP}
\left\{\begin{split}
&\frac{dU}{dt} = B U + F,\\
&U(0)=U^0,
\end{split}\right.
\end{equation}
where
\begin{enumerate}
  \item [$a)$] $U=[U_{ij}]$ is the solution vector:

  \noindent $U=(u_{22}, u_{23}, \ldots u_{2(N_y-1)}, u_{32}, u_{33}, \ldots u_{3(N_y-1)},\ldots u_{(N_y-1)2}, u_{(N_y-1)3}, \ldots u_{(N_y-1)(N_y-1)})$, and $U^0$ represents the initial solution vector.
  \item [$b)$] $F=[F_{ij}]$ is the vector of order $(N_y-2)(N_y-2)$ containing the boundary values, i.e.,

\noindent   $F_{ij}=\alpha_x  \(a_{i1}^{(2)} u_{1j}+  u_{N_x j}\)+ \alpha_y \(b_{j1}^{(2)} u_{i1}+b_{jN_y}^{(2)} u_{iN_y}\) - \beta_x \(a_{i1}^{(1)} u_{1j}+a_{iN_x}^{(1)} u_{N_xj}\) - \beta_y \(b_{j1}^{(1)} u_{i1}+ b_{j N_y}^{(1)} u_{iN_y}\).$

\item [$c)$] $B$ be a square matrix of order $(N_x-2)(N_y-2)$ given as
\begin{equation}\label{eq-B}
B = \alpha_x A_2 + \alpha_y B_2 -\beta_x A_1 - \beta_y B_1,
\end{equation}
where $A_{r}$ and $B_{r}$ $(r=1,2)$ are the block diagonal matrices of the weighting coefficients $a_{ij}^{(r)}$ and $b_{ij}^{(r)}$, respectively as given below
\begin{equation}\label{eq-ar-br}
\begin{array}{cc}
 A_r = \left[
    \begin{array}{cccc}
      a_{22}^{(r)} I & a_{23}^{(r)} I & \ldots & a_{2(N_x-1)}^{(r)} I \\
      a_{32}^{(r)} I & a_{33}^{(r)} I &\ldots  & a_{3(N_x-1)}^{(r)} I \\
      \vdots & \vdots & \ddots & \vdots \\
      a_{(N_x-1)2}^{(r)} I & a_{(N_x-2)3}^{(r)} I & \ldots  & a_{(N_x-1)(N_x-1)}^{(r)} I \\
    \end{array}
  \right], \mbox{ and }
   &  ~ B_r = \left[
    \begin{array}{cccc}
      M_r & O & \ldots & O \\
      O & M_r &\ldots  & O \\
      \vdots & \vdots & \ddots & \vdots \\
      O & O & \ldots  & M_r \\
    \end{array}
  \right],
\end{array}
\end{equation}
where $I$ and $O$ are the matrices of order $N_y-2$, and the sub-matrix $M_r$ of the block diagonal matrix $B_r$ is given by
\begin{equation}
\begin{array}{cc}
M_r = \left[
    \begin{array}{cccc}
      b_{22}^{(r)} & b_{23}^{(r)} & \ldots & b_{2(N_y-1)}^{(r)} \\
      b_{32}^{(r)} & b_{33}^{(r)} & \ldots & b_{3(N_y-1)}^{(r)} \\
      \vdots & \vdots & \ddots & \vdots \\
      b_{(N_y-1)2}^{(r)} & b_{(N_y-1)3}^{(r)} & \ldots & b_{(N_y-1)(N_y-1)}^{(r)}
    \end{array}
  \right].
\end{array}
\end{equation}
\end{enumerate}

Finally, we adopted SSP-RK54 scheme \cite{SR02} to solve the initial value system \eqref{eq-ode-IVP} as:
\begin{equation*}
\begin{split}
&U^{(1)} = U^m+ 0.391752226571890 \triangle t L(U^m) \\
&U^{(2)} = 0.444370493651235 U^m+ 0.555629506348765 U^{(1)} + 0.368410593050371 \triangle t L(U^{(1)}) \\
&U^{(3)}= 0.620101851488403 U^m+ 0.379898148511597 U^{(2)} + 0.251891774271694 \triangle t L(U^{(2)}) \\
&U^{(4)}= 0.178079954393132 U^m+ 0.821920045606868 U^{(3)} + 0.544974750228521 \triangle t L(U^{(3)}) \\
&U^{m + 1}= 0.517231671970585 U^{(2)}+ 0.096059710526147 U^{(3)} + 0.063692468666290 \triangle t L(U^{(3)})\\
& \qquad \quad + 0.386708617503269 U^{(4)} + 0.226007483236906 \triangle t L(U^{(4)}),
\end{split}
\end{equation*}
where $ L U = B U + F.$

\section{Stability of the methods for $2$D convection-diffusion equation} \label{sec-stab}
The stability of the method MTB-DQM for $2$D convection-diffusion equation \eqref{eq-CDEs} depends on the stability of
the initial value system of  ODEs as defined in \eqref{eq-ode-IVP}. Noticed that whenever the system of ODEs \eqref{eq-ode-IVP} is unstable, the proposed method for temporal discretization may not converge to the exact solution. Moreover, being the exact solution can directly obtained by means of the eigenvalues method, the stability of \eqref{eq-ode-IVP} depends on the eigenvalues of the coefficient matrix $B$ \cite{Jain83}.  In fact, the stability region is the set $\mathcal{S}= \{z \in C: \mid R(z)\mid \leq 1, z = \lambda_B \triangle t \}$, where $R(.)$ is the stability function and $\lambda_B$ is the eigenvalue of the coefficient matrix $B$. The stability region for SSP-RK54 scheme is depicted in \cite[Fig.1]{SPK16}, from which one can clam that for the stability of the system \eqref{eq-ode-IVP} it is sufficient that $\lambda_B \triangle t \in \mathcal{S}$ for each eigenvalue $\lambda_B$ of the coefficient matrix B. Hence, the real part of each eigenvalue is necessarily either zero or negative.

\begin{figure}
\centering
\includegraphics[height=4.5cm, width=6.25cm]{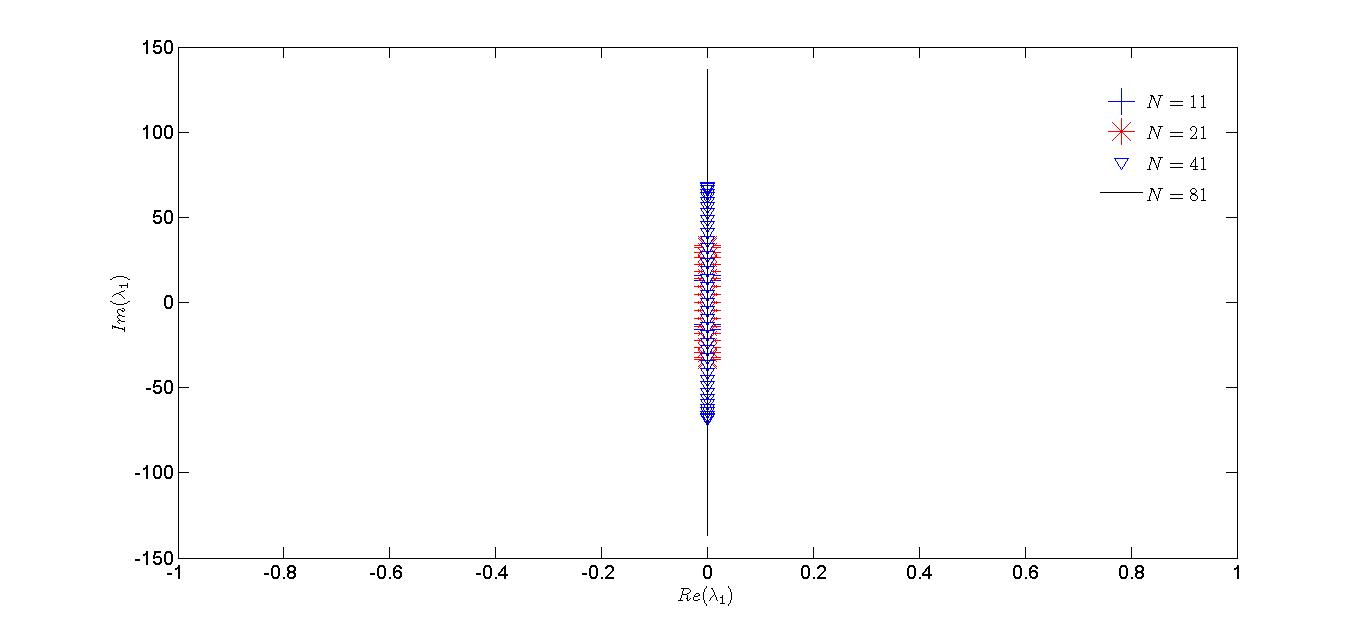}
\includegraphics[height=4.5cm, width=10.05cm]{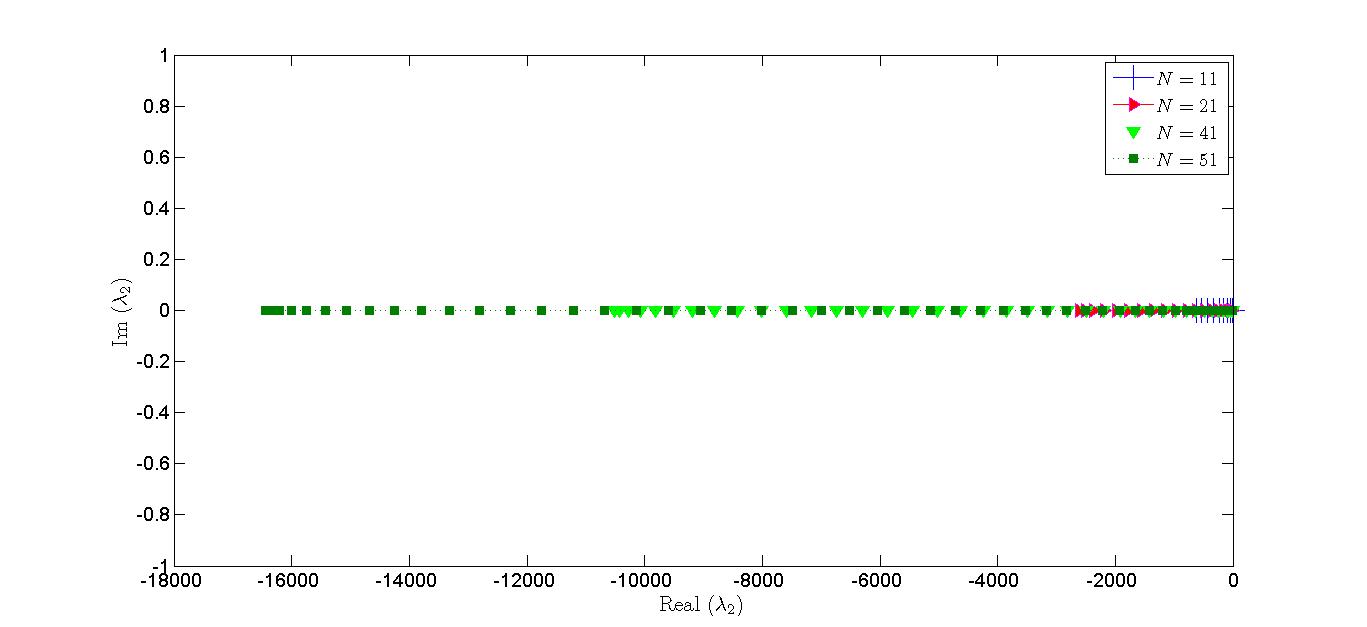}
\caption{Eigenvalues $\lambda_1$ and $\lambda_2$ for different grid sizes}\label{eq-eignv}

%\begin{figure}
%\centering
%\includegraphics[height=5.5cm, width=7.5cm]{A2x11}
%\includegraphics[height=5.5cm, width=7.5cm]{B2y11}
%\includegraphics[height=5.5cm, width=7.5cm]{A2x21}
%\includegraphics[height=5.5cm, width=7.5cm]{B2y21}
%\includegraphics[height=5.5cm, width=7.5cm]{A2x31}
%\includegraphics[height=5.5cm, width=7.5cm]{B2y31}
%\includegraphics[height=5.5cm, width=7.5cm]{A2x41}
%\includegraphics[height=5.5cm, width=7.5cm]{B2y41}
%\caption{Eigen values of $A_2$ (left) and $B_2$ (right) for different grid sizes}\label{eq-eigd2}
%\end{figure}
\includegraphics[height=5.15cm, width=15.0cm]{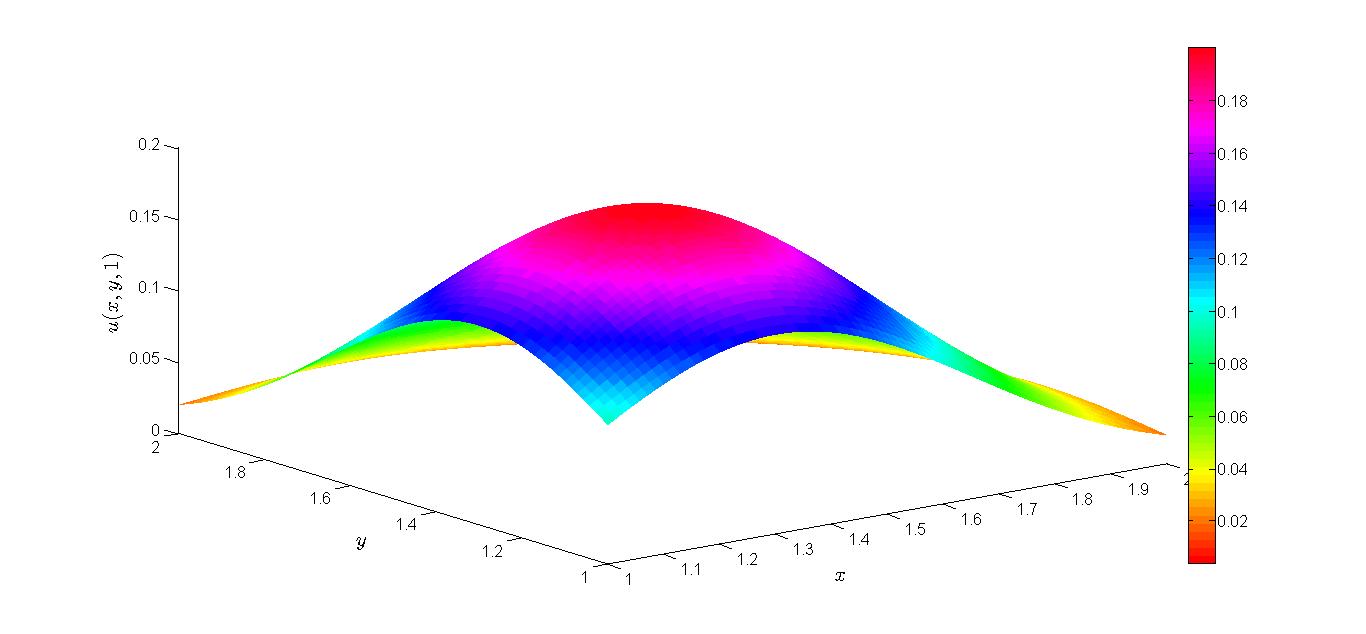}
\caption{Surface plot of distribution of $2$D convection-diffusion equation \eqref{eq-CDEs} with $\alpha_x=\alpha_y=0.05$, $\beta_x=\beta_y=0.8$ at $t = 1$}\label{ex1-fig1.1}
\includegraphics[height=5.15cm, width=15.0cm]{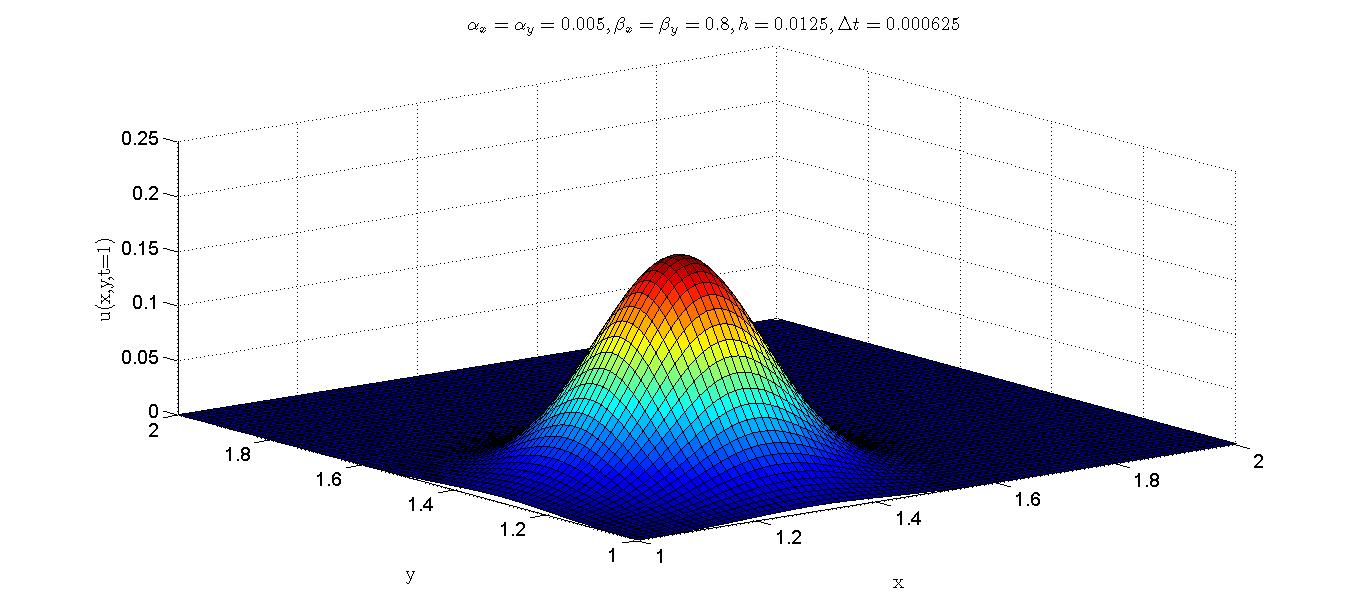}
\caption{Surface plot of distribution of $2$D convection-diffusion equation \eqref{eq-CDEs} with $\alpha_x=\alpha_y=0.005$, $\beta_x=\beta_y=0.8$ at $t = 1$ }\label{ex1-fig1.2}
\includegraphics[height=5.15cm, width=7.5cm]{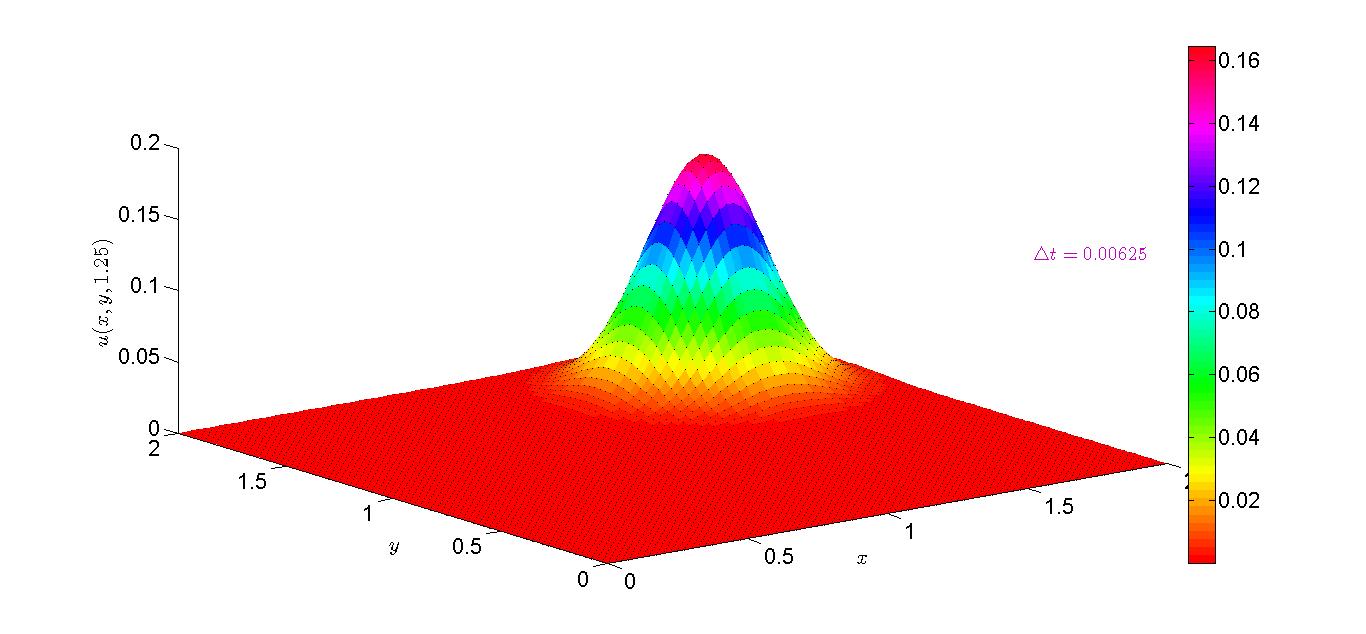}
\includegraphics[height=5.15cm, width=7.5cm]{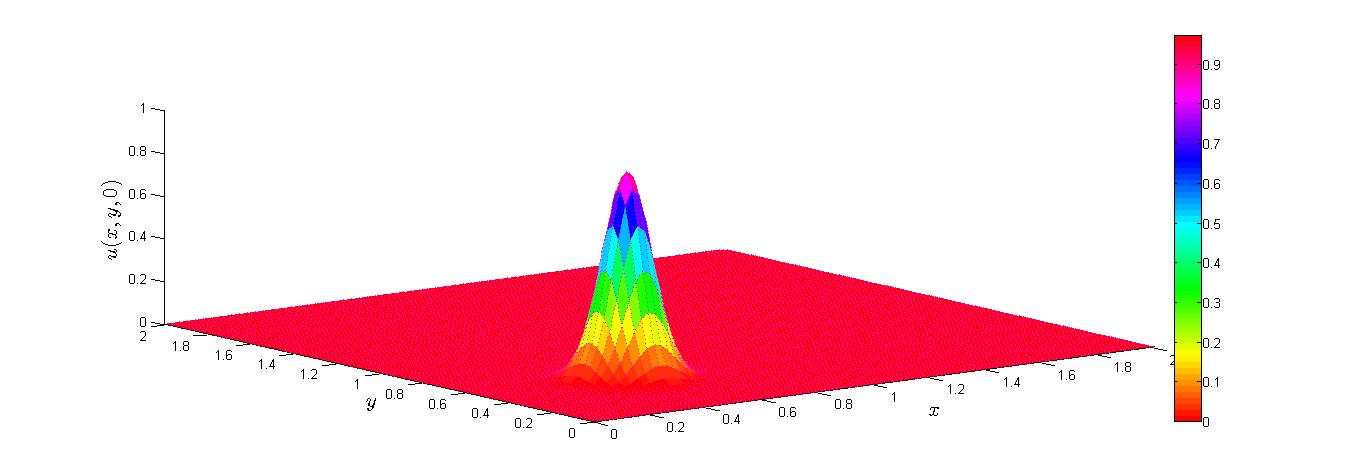}
\caption{ Surface plot of the inititial distribution (right) and the distribution at $t = 1.25$ (left) of 2D convection-diffusion equation \eqref{eq-CDEs} with $\alpha_x=\alpha_y=0.01$, $\beta_x=\beta_y=0.8$ for $h_x=h_y=0.025$ }\label{ex1-fig1.3}
\end{figure}

It is seen that the eigenvalues of the matrices $A_r$ and $B_r$ $(r=1, 2)$ have identical nature. Therefore, it is sufficient to compute the eigenvalues, $\lambda_1$ and $\lambda_2$, of the matrices  $A_1$ and $A_2$ for different values of grid sizes $(N_x \times N_y)$. The eigenvalues $\lambda_1$ and $\lambda_2$ for $N_x = N_y=11, 21, 31, 41$ has been depicted in Figure \ref{eq-eignv}. Analogously, one can compute the eigenvalues $\lambda_1$ and $\lambda_2$ using mECDQ method \cite{SPK16} or mExp-DQM. It is seen that in either case, $\lambda_1$ and $\lambda_2$ have same nature as in Figure \ref{eq-eignv}.  Further, we can get from Figure \ref{eq-eignv} that each eigenvalue $\lambda_B$ of the matrix $B$ as defined in Eq. \eqref{eq-B} is real and negative. This confirms that the proposed methods produces stable solutions for two dimensional convection-diffusion equations.

\section{Numerical results and discussions} \label{sec-num}
This section deals with the main goal, the numerical study of three test problems of the initial value system of convection-diffusion equations with both kinds of the boundary conditions has been done by adopting the methods MTB-DQM, mExp-DQM and mECDQ method along with the integration SSP-RK54 scheme. The accuracy and the efficiency of the methods have been measured in terms of the discrete error norms: namely- average $|\mbox{error}|$ norm ($L^2$-error norm) and the maximum error ($L_{\infty}$ error norm). %In this section, we take $h_x=h_y=h$.

\begin{problem}\label{ex1}
Consider the initial value system of $2$D convection-diffusion equation \eqref{eq-CDEs} with $u_0(x, y) = \exp\left\{-\frac{(x-x_0)^2}{\alpha_x} -\frac{(y-y_0)^2}{\alpha_y} \right\}$, while values of $f_i$ for $1\le i \le 4$ can be extracted from the exact solution
\begin{equation*}
u(x, y, t) =\frac{1}{{1+4t}} \exp\left\{- \frac{(x-x_0 - \beta_x t)^2}{\alpha_x(1+4t)} - \frac{(y-y_0 - \beta_y t)^2}{\alpha_y(1+4t)} \right\},
\end{equation*}
where initial condition is a Gaussian pulse with unit hight centered at $(x_0, y_0)=(0.5, 0.5)$.

We have computed the numerical solution of Problem \ref{ex1} for $\beta_x = \beta_y = 0.8$ and different values of $\alpha_x, \alpha_y$.

For $\Omega = [1, 2]^2$: $L_\infty$ errors and the rate of convergence (ROC) \cite{SPK16} at $t = 1$ for different values of $h^2=\triangle t$ ($h_x=h_y=h$) has been reported in Table \ref{tab1.1} for $\alpha_x = \alpha_y=0.05, 0.005$. Table  \ref{tab1.1} confirms that the proposed solutions obtained by these methods are more accurate, and approaching towards exact solutions.  The behavior of the approximate solution at time $t = 1$ taking $h = 0.0125, \triangle t = 0.000625$ is depicted in Figure \ref{ex1-fig1.1} for $\alpha_x=\alpha_y=0.05$ and in Figure \ref{ex1-fig1.2} for $\alpha_x=\alpha_y=0.005$.

For $\Omega =[0, 2]^2$, time $t=1.25$, and grid space step size $h = 0.025$: the $L^2$ and $L_\infty$ error norms in the proposed solutions have been compared with the errors in the solutions by various schemes in \cite{TG07,NT88,KDD02,KZ04,TG07,DM08} in Table \ref{tab1.2} for $\triangle t = 0.00625$ and in Table \ref{tab1.3} for $\triangle t = 0.0125$. The initial solution and MTB-DQM solutions in $[0, 2]^2$ with the parameter values $\alpha_x = \alpha_y= 0.01, t = 1.25$ have been depicted in Figure \ref{ex1-fig1.3}, and similar behavior is seen from the other two methods. It is evident from the above reports that the proposed results are more accurate as compared to \cite{TG07,NT88,KDD02,KZ04,TG07,DM08} and approaching towards the exact solutions.
\begin{table}[!b]
\caption{Rate of convergence of $L_\infty$ error norms in Problem \ref{ex1} in $[1,2]^2$ with $\alpha_x=\alpha_y=0.05, 0.005, \beta_x=\beta_y=0.8$ taking $\triangle t=h^2$ and various grid sizes $N_x\times N_y$ with $N_x= N_y$}  \label{tab1.1}
\vspace{.2cm}
\centering
\begin{tabular}{l*8l}
\toprule
\multicolumn{9}{l}{$\alpha_x=\alpha_y=0.005 $} \\
\hline
 $N_x$  & \multicolumn{2}{l}{MTB-DQM} & {} & \multicolumn{2}{l}{mExp-DQM (p=0.001) } & {} & \multicolumn{2}{l}{mECDQ $(\lambda = -0.30)$ } \\ \cline{2-3} \cline{5-6} \cline{8-9}	
  & $L_\infty$  $\quad$& ROC $\quad$  & {} &	$L_\infty$  $\quad$& ROC $\quad$  & {} & $L_\infty$  $\quad$& ROC\\
\midrule
5 &2.91E-03  &	   & {}  &2.91E-03&       & {} &1.94E-03    &   	\\
10&6.41E-04  &2.18 & {}  &6.42E-04&2.181  & {} &4.93E-04	&1.98\\
20&6.50E-05  &3.30 & {}  &6.51E-05&3.303  & {} &6.05E-05	&3.03\\
40&7.19E-06  &3.18 & {}  &7.18E-06&3.180  & {} &6.91E-06	&3.13\\
80$\qquad \qquad$ &7.22E-07  $\qquad \qquad $ &3.31 &$\qquad$  &7.00E-07$\qquad \qquad $&3.358  & $\qquad$ &6.99E-07	$\qquad $&3.31\\
\multicolumn{9}{l}{$\alpha_x=\alpha_y=0.05 $} \\
\midrule
5 &4.93E-04&		       &&4.97E-04 &       & {} &3.20E-04&	 \\
10&4.86E-05&	3.34	   &&4.87E-05 &3.35	  & {} &4.39E-05&2.87\\
20&4.70E-06&	3.37	   &&4.70E-06 &3.37   & {} &4.35E-06&3.34\\
40&4.51E-07&	3.38	   &&4.46E-07 &3.40   & {} &4.50E-07&3.27\\
80$\qquad \qquad$&4.58E-08$\qquad \qquad$ & 3.30	  $\qquad \qquad$  &&3.49E-08 $\qquad \qquad$&3.67	  & {} &4.37E-08$\qquad \qquad$ &3.36\\
\bottomrule
\end{tabular}
\end{table}

\begin{table}[!b]
\caption{Comparison  $L^2$ and $L_\infty$ error norms for Problem \ref{ex1} with $\alpha_x=\alpha_y=0.01; \beta_x=\beta_y=0.8$ for $h_x=h_y=0.025, \triangle t = 0.00625$}  \label{tab1.2}
\centering
\begin{tabular}{*3l}
\toprule
Schemes		&		$L^2$			&	$L_\infty$ \\%\cline{2-4}\cline{6-8} \cline{10-12}\cline{14-16}
\midrule
P-R-ADI \cite{TG07}				& 3.109E-04		$\qquad\qquad $	 $\qquad\qquad $	 &7.778E-03	\\
Noye $\&$ Tang \cite{NT88}				& 1.971E-05				&6.509E-04	 \\
Kalila et al. \cite{KDD02}			& 1.597E-05				&4.447E-04	\\
Kara $\&$ Zhang ADI \cite{KZ04}		& 9.218E-06				&2.500E-04	\\
Tang $\&$ Ge ADI\cite{TG07}			& 9.663E-06				&2.664E-04	\\
Dehghan $\&$ Mohebbi \cite{DM08}$\qquad\qquad $  $\qquad\qquad $	& 9.431E-06				 &2.477E-04	\\
MTB-DQM					        & 8.026E-12			    &4.327E-08	\\
mExp-DQM(p=0.0001)					        & 7.030E-12			    &4.154E-08	 \\
mECDQ ($\lambda=-0.004$)					        & 4.504E-12			    &3.343E-08	 \\
MCB-DQM	($\lambda=0$)				        & 8.269E-12			    &4.388E-08	 \\
\bottomrule
\end{tabular}
\end{table}
\begin{table}[!b]
\caption{Comparison  $L^2$ and $L_\infty$ error norms for Problem \ref{ex1} with $\alpha_x=\alpha_y=0.01; \beta_x=\beta_y=0.8$ for $h_x=h_y=0.025, \triangle t= 0.0125$}  \label{tab1.3}
\centering
\begin{tabular}{*3l}
\toprule
Schemes		&		$L^2$			&	$L_\infty$ \\
\midrule
Noye $\&$ Tang \cite{NT88}				&1.430E-05	&4.840E-04\\
Kalila et al.\cite{KDD02}				&1.590E-05	&4.480E-04\\
Dehghan $\&$ Mohebbi \cite{DM08} $\qquad\qquad $  $\qquad\qquad $ &	 9.480E-06	 $\qquad\qquad $  $\qquad\qquad $ $\qquad\qquad $ &			 2.469E-04\\
MTB-DQM	& 8.900E-12			& 4.501E-08 \\
mExp-DQM(p=0.0001)	& 9.147E-12			& 4.315E-08 \\
mECDQ ($\lambda=-0.005$)	& 4.646E-12			& 3.852E-08 \\
MCB-DQM	($\lambda=0$) & 9.147E-12			& 4.561E-08 \\
\bottomrule	
\end{tabular}
\end{table}
\end{problem}

\begin{figure}
%\centering
%\includegraphics[height=5.5cm, width=7.5cm]{fig36cont} ?????
%\includegraphics[height=5.5cm, width=7.5cm]{fig36surf}
%\caption{ Contour plot (left) and surface plot (right) of the distribution of 2D convection-diffusion equation \eqref{eq-CDEs} with $\alpha_x=\alpha_y=0.10$, $\beta_x=\beta_y=1$  at $t = 1$ }\label{ex2-fig2.1}
\includegraphics[height=5.15cm, width=7.5cm]{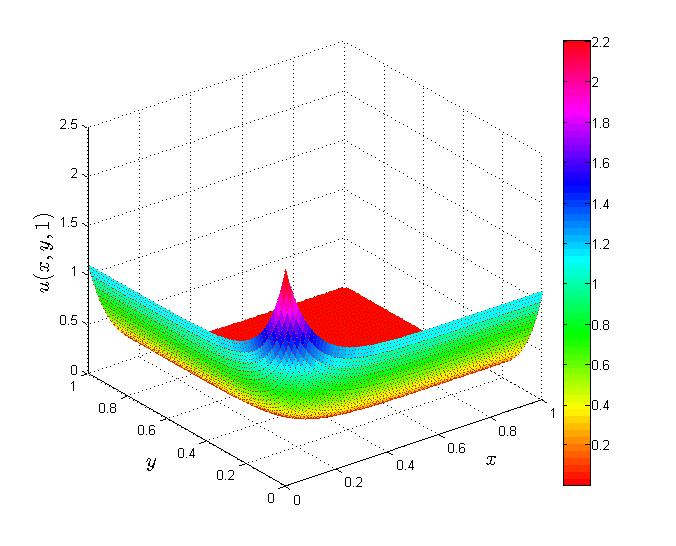}
\includegraphics[height=5.15cm, width=7.5cm]{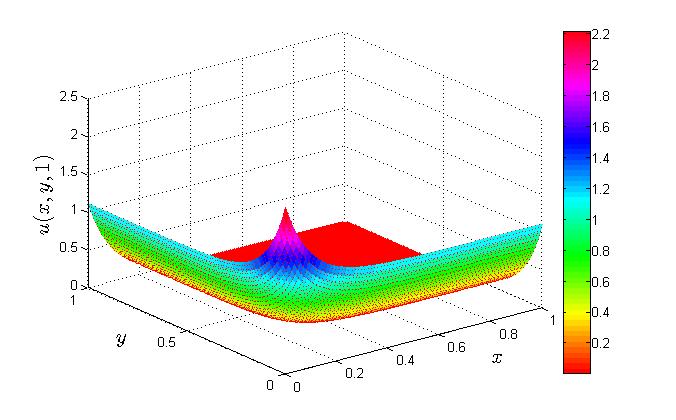}
\caption{The approximate (left) and exact solution behavior (right) of  $2$D convection-diffusion equation \eqref{eq-CDEs} with $\alpha_x=\alpha_y=0.10$, $\beta_x=\beta_y=-1$  at $t = 1$ }\label{ex2-fig2.2}

\includegraphics[height=5.15cm, width=7.5cm]{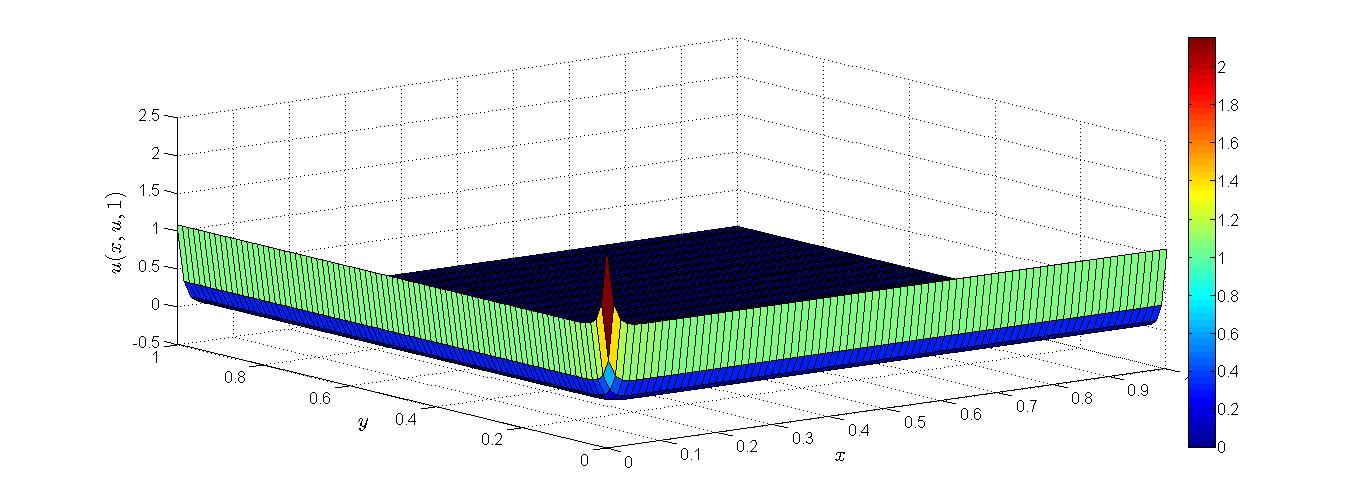}
\includegraphics[height=5.15cm, width=7.5cm]{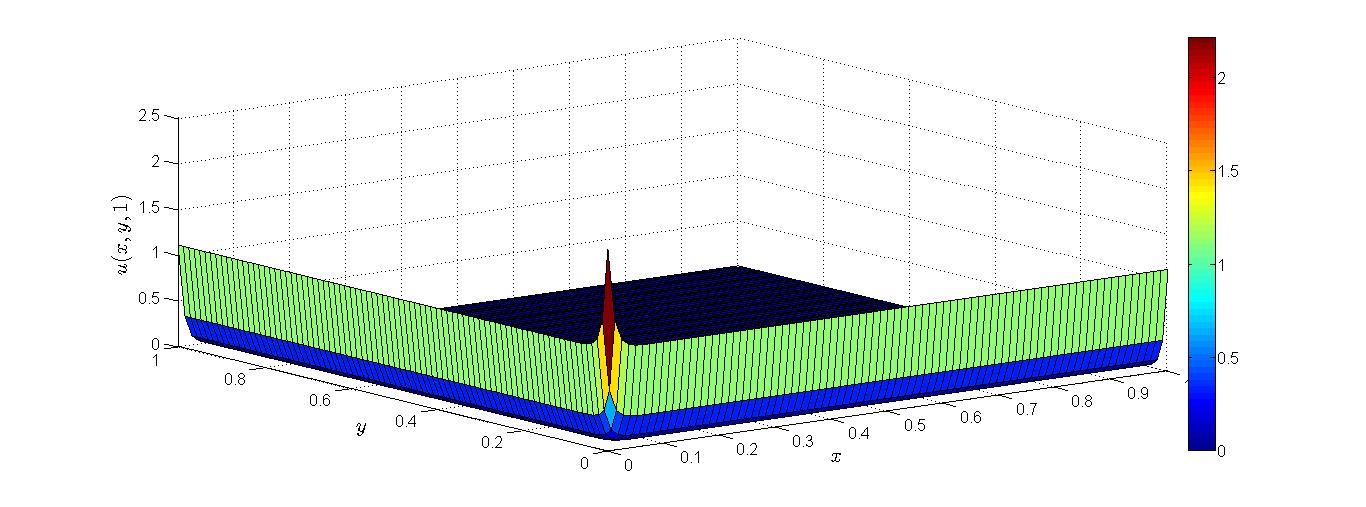}
\caption{The approximate solution behavior of $2$D convection-diffusion equation \eqref{eq-CDEs} with $\alpha_x=\alpha_y=0.01$, $\beta_x=\beta_y=-1$  at $t = 1$ }\label{ex2-fig2.3}

%\includegraphics[height=5.15cm, width=7.5cm]{fig3p91}
%\includegraphics[height=5.15cm, width=6.5cm]{fig3p92}
%\caption{Contour plot of the initial distribution of  2D convection-diffusion equation \eqref{eq-CDEs} in Example \ref{ex3}}\label{ex3-fig3.1}

\includegraphics[height=5.15cm, width=7.75cm]{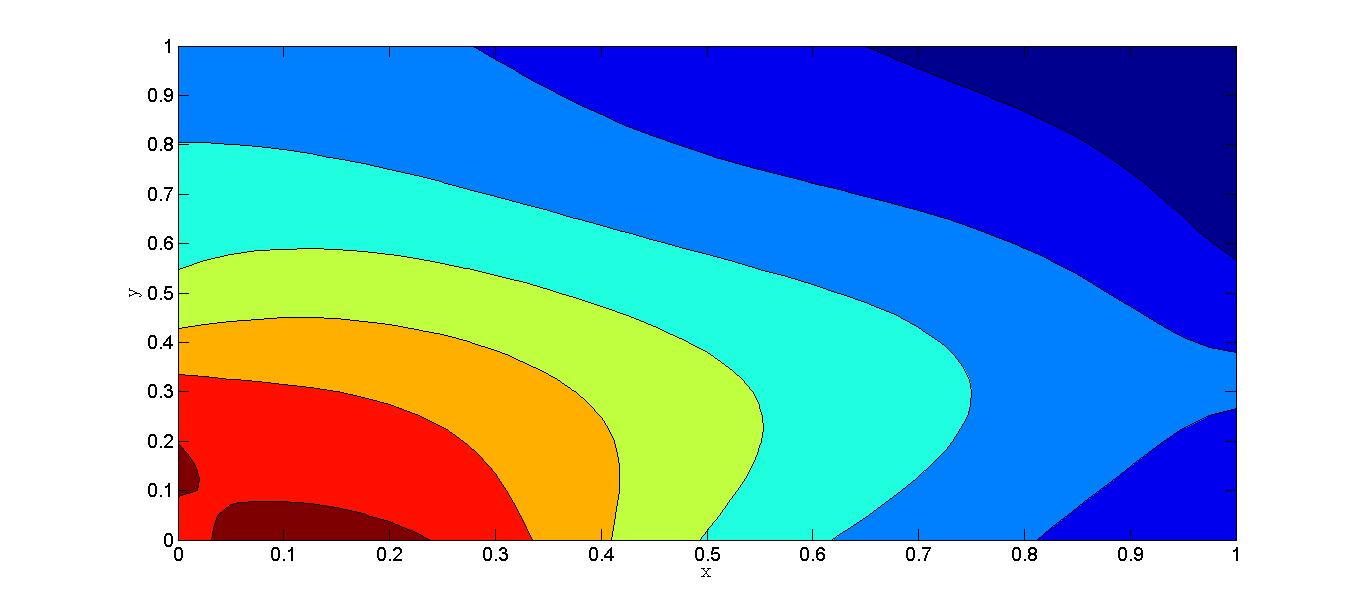}
\includegraphics[height=5.15cm, width=7.75cm]{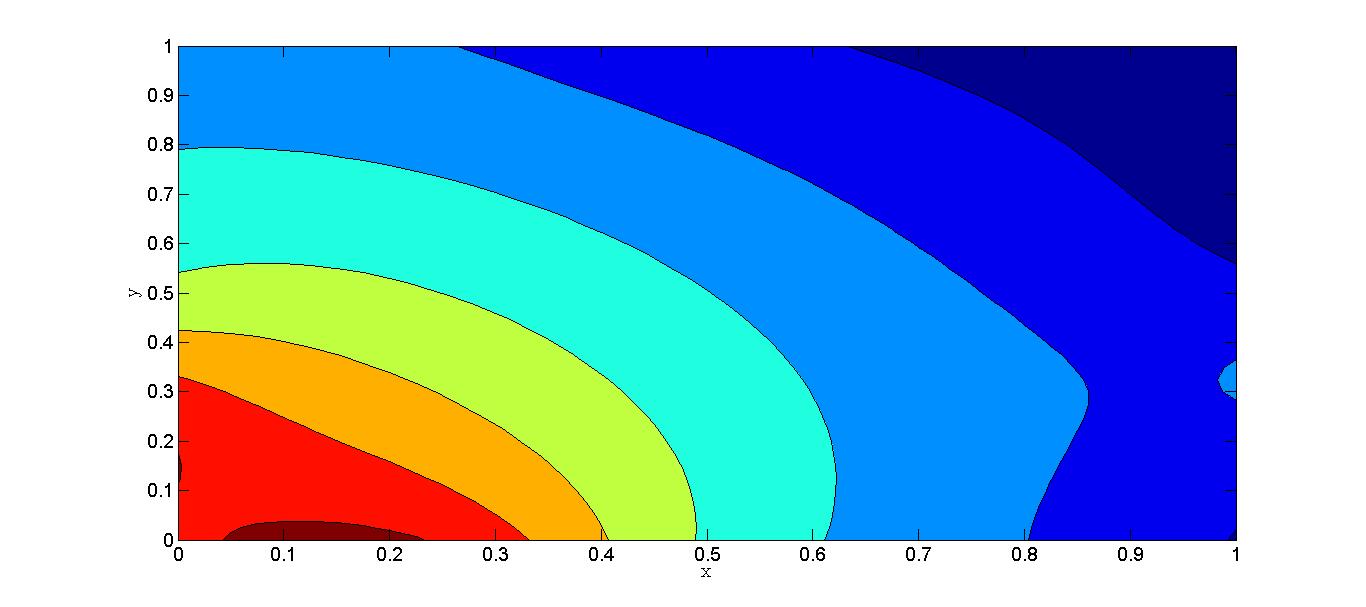}
\caption{The approximate solution behavior of the distribution of $2$D convection-diffusion equation \eqref{eq-CDEs} in Example \ref{ex3} at $t = 0.1$ (left) and $ t=0.5$ (right)} \label{ex3-fig3.2}
\end{figure}

\begin{problem}\label{ex2}
Consider the initial value system \eqref{eq-CDEs} of $2$D convection-diffusion equation with $\Omega= [0, 1]^2$ with $u_0(x, y) = a\( \exp(-c_x x) + \exp(-c_y y)\) $, where
\begin{equation*}
c_x =\frac{ -\beta_x\pm \sqrt{\beta_x^2+4b\alpha_x}}{2 \alpha_x}>0, \mbox{ and } c_y =\frac{ -\beta_y\pm \sqrt{\beta_y^2+4b\alpha_y}}{2 \alpha_y}>0,
\end{equation*} and the Neumann boundary condition
\begin{equation*}
\left\{ \begin{split}
&\left.\frac{\partial u }{\partial x}\right|_{x=0, y} = -a c_x \exp(bt), \quad \left.\frac{\partial u }{\partial x}\right|_{x=2, y} = -ac_x \exp(bt-c_x),\\
&\left.\frac{\partial u }{\partial y}\right|_{x, y=0} = -a c_y \exp(bt), \quad \left.\frac{\partial u }{\partial y}\right|_{x, y=2} = -a c_y \exp(bt-c_y),
\end{split}
\right. \qquad \quad (x, y, t)\in \partial\Omega \times (0, T],
\end{equation*}or the Dirichlet's conditions can be extracted from the exact solution \cite{CDN06}:
\begin{equation*}
u(x, y, t)= a \exp(bt)\( \exp(-c_x x) + \exp(-c_y y)\).
\end{equation*}

The computation by the proposed methods has been done for different values of $\alpha_x=\alpha_y$ and $\beta_x=\beta_y$.

For $\beta_x=\beta_y=1$, we take $\triangle t=0.001$ for the solutions at $t=1$. The rate of convergence and $L^2$ error norms in the proposed solutions has been compared with that of due to fourth-order compact finite difference scheme \cite{DM08} for $\alpha_x=\alpha_y = 0.1, 0.01$, in Table \ref{tab2.1}.  It is found that the proposed solutions from either method are more accurate in comparison to \cite{DM08}, and are in good agreement with the exact solutions. %The surface plots and contour plots of MTB-DQM solutions at $t=1$, and $h=0.025$ have been depicted in Figure \ref{ex2-fig2.1}.

For $\beta_x=\beta_y=-1$, we take $\triangle t=0.0005$ for the solution at $t=1$. In Table \ref{tab2.2}, the $L^2$ and $L_\infty$ error norms are compared with that obtained by fourth-order compact finite difference scheme \cite{DM08} for $\beta_x=\beta_y=-1$, and $\alpha_x=\alpha_y=0.01$. Rate of convergence have been reported in Table \ref{tab2.3}. From Tables \ref{tab2.2} and \ref{tab2.3}, we found that the proposed solutions are more accurate as compared to the results in \cite{DM08}, while the rate of convergence is linear. The behavior of solutions is depicted in Figure \ref{ex2-fig2.2} and \ref{ex2-fig2.3} with $h=0.025$ for $\alpha_x=\alpha_y=0.1$ and $\alpha_x=\alpha_y= 0.01$, respectively.

For the same parameter, mentioned above, the numerical solution is obtained by using Neumann conditions and reported in Table \ref{tab2.4}. The obtained results are in good agreement with the exact solutions, the rate of convergence for $L\infty$ is quadratic for each method.

\begin{table}[!b]
\caption{$L^2$-error norms and ROC for Problem \ref{ex2} with $\beta_x=\beta_y=1$ and $\alpha_x=\alpha_y=0.10, 0.01$ at $t=1$ } \label{tab2.1}
\centering
\begin{tabular}{*2lclcllclcllclcllc}
\toprule
 \multicolumn{18}{l}{ $\alpha_x=\alpha_y = 0.10$ and $\triangle t =0.0005$  }  \\
\hline
{} & \multicolumn{4}{l}{MTB-DQM} &&  \multicolumn{4}{l}{mExp-DQM $(p=1)$} && \multicolumn{4}{l}{ mECDQ $(\lambda=-0.0001)$} &&  \multicolumn{2}{l}{CFDS4\cite{DM08}}  \\ \cline{2-5}\cline{7-10}\cline{12-15}\cline{17-18}
$h$& $L\infty$ &	ROC &	$L^2$ &ROC&& 	$L\infty$  &	ROC & $L^2$ &	 ROC &&	$L\infty$  &ROC& 	$L^2$  &	ROC && 	$L^2$  &	ROC\\
\midrule 	
0.2	  & 4.231E-06& 	    &   4.969E-10&	    &&	1.714E-05&		&  3.584E-09&		 &&   4.231E-06&       &  4.999E-10  &   &&   &  \\
0.1	  & 5.155E-07&	3.0&	2.472E-11&	4.3&&	5.168E-07&	5.1&  2.473E-11&	 7.2&&	5.153E-07&	3.0&	2.478E-11&	4.3	&&2.289E-10&	\\
0.05  &	5.108E-08&	3.3&	9.041E-13&	4.8&&	5.105E-08&	3.3&	 9.044E-13&	 4.8&&	5.107E-08&	3.3&	9.049E-13&	4.8&&1.621E-11&	3.82\\
0.025 &	1.147E-08&	2.2&	1.922E-14&	5.6&&	1.147E-08&	2.2&	 1.922E-14&	 5.6&&	1.147E--08&	2.4&	1.923E-14&	5.6&&8.652E-13&	4.23\\
\multicolumn{18}{l}{  $\alpha_x=\alpha_y = 0.01,  \triangle t =0.001$ } \\
\midrule 	
0.2	 & 7.4031E-06&		&  1.4907E-09&		&&  3.5634E-05&		&  1.5314E-08&		 &&  2.0585E-06&		&   6.7290E-11&		&&		      &  \\
0.1	 & 1.1460E-06&	2.7	&  1.3103E-10&	3.5	&&  1.1421E-06&	5.0	&  1.3148E-10&	 6.9	&&  1.5919E-07&	3.7	&   1.1624E-12&	5.9	&&	2.749E-09&	  \\
0.05 & 1.4578E-07&	3.0	&  7.5377E-12&	4.1	&&  1.4582E-07&	3.0	&  7.5396E-12&	 4.1	&&  1.9012E-08&	3.1	&   5.4939E-14&	4.4	&&	2.394E-10&	 3.52\\
0.025& 1.3902E-08&	3.4	&  2.7129E-13&	4.8	&&  1.3904E-08&	3.4	&  2.7132E-13&	 4.8	&&  6.2220E-09&	1.6	&   9.0456E-15&	2.6	&&	1.658E-11&	 3.85\\
\bottomrule
\end{tabular}
\end{table}

\begin{table}[!b]
\caption{Comparison of the error norms in Problem \ref{ex2} with $\beta_x=\beta_y=-1$, $\alpha_x=\alpha_y=0.01 $ taking $\triangle t=0.0005$ with CFDS4\cite{DM08}} \label{tab2.2}
\centering
\begin{tabular}{*6l*6l*6l}
\toprule
{} & \multicolumn{3}{l}{MTB-DQM}&{}&\multicolumn{3}{l}{mExp-DQM,$p=0.0001$}&{}&\multicolumn{3}{l}{mECDQ,$\lambda=-1.75$}& {} & \multicolumn{3}{l}{CFDS4 \cite{DM08}}\\ \cline{2-4}\cline{6-8}\cline{10-12}\cline{14-16}
$h$  &  $L_\infty$ &{} &	$L^2$ & {} & $L_\infty$ &{} &	$L^2$ & {} & $L_\infty$ &{} &	$L^2$ & {} & $L_\infty$ &{} &	$L^2$ \\
\midrule 							
0.04 $\quad$ & 	4.1718E-03$\quad$ &&	1.7833E-03$\quad$ &&		 1.7475E-03$\quad$ &&7.9168E-04&&		1.2939E-03$\quad$ &&	3.6957E-04&&		 1.1826E-01	$\quad$ &&	4.9331E-03\\			
0.02 &	5.2095E-04&&	6.8410E-05&&		5.2097E-04&&	6.8445E-05&&		 6.3735E-04&&	 6.8622E-05&&		1.5310E-02&&	4.1351E-04\\		
0.01 &	3.2496E-04&&	2.1508E-05&&		 3.2484E-04	&&	 2.1493E-05&&		 2.5394E-04&&	1.8260E-05&&		9.4696E-04&&	2.8405E-05 \\				
\bottomrule
\end{tabular}
\end{table}

\begin{table}[!b]
\caption{$L^2$ and $L\infty$ error norms in Problem \ref{ex2} with $\beta_x=\beta_y=-1$ and different values of $\alpha_x=\alpha_y$} \label{tab2.3}
\centering
\begin{tabular}{*5l*5l*5l}
\toprule
{} & \multicolumn{8}{l}{$\alpha_x=\alpha_y = 0.01, \triangle t =0.0005$} \\\cline{2-15}
{} & \multicolumn{4}{l}{MTB-DQM} &&  \multicolumn{4}{l}{mExp-DQM,$p=0.0001$}& &\multicolumn{4}{l}{mECDQ, $\lambda=-0.9$}  \\ \cline{2-5}\cline{7-10}\cline{12-15}
$h$   & $L\infty$ &	ROC &	$L^2$      &  ROC  && 	$L\infty$  &	ROC & $L^2$ &	 ROC &&	$L\infty$  &ROC& 	$L^2$  &	ROC\\
\midrule 	
0.1  &	2.6263E-02&		&  3.6277E-02	&	      &&2.6280E-02&		&   3.6670E-02&		  &&2.4464E-02&		&4.0175E-02	&	\\
0.05&	6.8144E-03	&1.9&	4.1447E-03	&3.1	&&6.8144E-03	&1.9&	 4.1339E-03&  3.1  &&3.3580E-03&	2.9 &1.6264E-03	&4.6	\\
0.025&	1.1634E-03	&2.6&	2.7780E-04	&3.9	&&1.0719E-03	&2.7&	 1.9697E-04&  4.4  &&8.1776E-04&	2.0 &1.0090E-04	&4.0	\\
{} & \multicolumn{8}{l}{$\alpha_x=\alpha_y = 0.1, \triangle t =0.0005$} \\\cline{2-15}
{} & \multicolumn{4}{l}{} &&  \multicolumn{4}{l}{$p=0.0001$}& &\multicolumn{4}{l}{$\lambda=0.3$}  \\ \cline{2-5}\cline{7-10}\cline{12-15}
0.01 &5.0998E-04&		&3.5191E-03&		&&5.0738E-04&		&1.0116E-05&		 &&7.4180E-04&		&2.5558E-05&	\\
0.05 &1.8776E-04&	1.4	&1.8193E-03&	1.0	&&1.8776E-04&	1.4	&2.9979E-06&	 1.8	&&2.0972E-04&	1.8	&5.1381E-06&	2.3\\
0.025 &9.9342E-05&	0.9	&9.1372E-04&	1.0	&&9.9342E-05&	0.9	&7.9462E-07&	 1.9	&&9.5677E-05&	1.1	&8.1188E-07&	2.7\\
\bottomrule
\end{tabular}
\end{table}

\begin{table}[!b]
\caption{ The $L^2$-error norm in Problem \ref{ex2} with Neumann boundary conditions for the different values of $\alpha_x, \alpha_y$	and $\beta_x, \beta_y$ with $ \triangle t =0.0005$} \label{tab2.4}
\begin{tabular}{l*4c*5c*5c}
\toprule
{} & \multicolumn{14}{c}{$\alpha_x=\alpha_y = 0.1$and $\beta_x=\beta_y=-1$} \\\cline{2-15}
{} & \multicolumn{4}{c}{MTB-DQM} &&  \multicolumn{4}{c}{mExp-DQM,$p=0.001$}& &\multicolumn{4}{c}{mECDQ, $\lambda=-0.9$}  \\ \cline{2-5}\cline{7-10}\cline{12-15}
$h$   & $L\infty$ &	ROC &	$L^2$      &  ROC  && 	$L\infty$  &	ROC & $L^2$ &	 ROC &&	$L\infty$  &ROC& 	$L^2$  &	ROC\\
\midrule 	
0.05	&1.0176E-02&		&6.7504E-03&		&&1.0177E-02&		 &6.7513E-03&		 &&9.1549E-03&		&5.7412E-03&	\\
0.025	&2.2164E-03&	2.2	&1.0428E-03&	2.7	&&2.2165E-03	&2.2 &1.0429E-03&	 2.7	&&1.9834E-03&	2.2	&8.7214E-04&	2.7\\
0.0125	&2.9460E-04&	2.9	&7.2888E-05&	3.8	&&2.9460E-04	&2.9 &7.2889E-05&	 3.8	&&2.5565E-04&3.0	&5.6518E-05&	3.9\\
%\bottomrule
{} & \multicolumn{14}{c}{$\alpha_x=\alpha_y = 0.01$ and $\beta_x=\beta_y=-1$} \\\cline{2-15}
{} & \multicolumn{4}{c}{} &&  \multicolumn{4}{c}{$p=0.001$}& &\multicolumn{4}{c}{ $\lambda=-0.9$}  \\ \cline{2-5}\cline{7-10}\cline{12-15}
0.05  &1.2070E-01&		&6.3285E-01&		&&1.2070E-01&		&6.3293E-01&		 &&1.0921E-01&		&6.0454E-01&	\\
0.025 &2.7808E-02&	2.1	&6.2302E-02&	3.3	&&2.7809E-02&	2.1	&6.2304E-02&	 3.3	&&2.6004E-02&	2.1	&5.5004E-02&	3.5\\
0.0125&6.0348E-03&	2.2	&6.4433E-03&	3.3	&&6.0348E-03&	2.2	&6.4433E-03&	 3.3	&&5.4733E-03&	2.2	&5.5668E-03&	3.3\\
{} & \multicolumn{14}{c}{$\alpha_x=\alpha_y = 0.1$ and $\beta_x=\beta_y=1$} \\\cline{2-15}
{} & \multicolumn{4}{c}{} &&  \multicolumn{4}{c}{p=0.01 }& &\multicolumn{4}{c}{$\lambda=-0.9$}  \\ \cline{2-5}\cline{7-10}\cline{12-15}	
0.05	&1.4931E-05&		&5.8044E-08&		&&1.4932E-05&		 &5.8056E-08&		 &&1.3598E-05&		&4.8332E-08&	\\
0.025	&3.3009E-06&	2.2	&1.0664E-08&	2.4	&&3.3009E-06&	2.2	 &1.0665E-08&	 2.4	&&2.9824E-06&	2.2	&8.7286E-09&	2.5\\
0.0125	&4.5947E-07&	2.8	&8.0780E-10&	3.7	&&4.5947E-07&	2.8	 &8.0780E-10&	 3.7	&&3.9840E-07&	2.9	&6.0761E-10&	3.8\\
{} & \multicolumn{14}{c}{$\alpha_x=\alpha_y = 0.01$ and $\beta_x=\beta_y=1$} \\\cline{2-15}
{} & \multicolumn{4}{c}{} &&  \multicolumn{4}{c}{p=0.01 }& &\multicolumn{4}{c}{$\lambda=-0.9$}  \\ \cline{2-5}\cline{7-10}\cline{12-15}	
0.05   &1.7132E-05&		    &6.9905E-08&		&&1.7132E-05&		 &6.9910E-08&		 &&1.6511E-05&		&6.5219E-08&	\\
0.025  &3.9931E-06&	   2.1 	&1.4135E-08&	2.3	&&3.9932E-06&	2.1	 &1.4136E-08&	 2.3	&&3.7892E-06&	2.1	&1.2756E-08&	2.4\\
0.0125 &8.8931E-07&	   2.2	&2.7052E-09&	2.4	&&8.8931E-07&	2.2	 &2.7052E-09&	 2.4	&&8.2511E-07&	2.2	&2.3316E-09&	2.5\\
\bottomrule
\end{tabular}
\end{table}

\end{problem}

\begin{problem}\label{ex3}
The initial value system of $2$D convection-diffusion equation \eqref{eq-CDEs} with $\Omega=[0, 1]^2$ and $u_0(x, y) = h(x, y)$, and
\begin{equation*}
\begin{split}
&f_1=h(0, y)-\overline{c}t;  \quad f_1=h(1, y)-\overline{c}t; \quad  f_3==h(x, 0)-\overline{c}t; \quad f_4=h(x, 1)-\overline{c}t,
\end{split}
\end{equation*}
where
\begin{equation*}
\begin{split}
h(x, y)&= 5 \exp\left\{-\frac{(9x-2)^2}{4} -\frac{(9y-2)^2}{4}\right\} + 7 \exp\left\{-\frac{(9x+1)^2}{50} -\frac{(9y+1)}{10}\right\}\\
&\qquad  + 4 \exp\left\{-\frac{(9x-7)^2}{4} -\frac{(9y-3)^2}{4}\right\} -2 \exp\left\{-(9x-4)^2 -(9y-7)^2\right\}
\end{split}
\end{equation*}
The distribution of the initial solution is depicted in Figure \ref{ex3-fig3.1}. The solutions behavior is obtained for the parameter values: $\alpha_x=0.2, \alpha_y=0.3, \beta_x=-0.1, \beta_y=0.2$ and $\triangle t=0.0005, h_x=h_y=0.025$, and is depicted in Figure \ref{ex3-fig3.2} due to MTB-DQM, also we noticed the similar characteristics obtained using mExp-DQM and mECDQ method. The obtained characteristics agreed well as obtained in \cite{DM08,ZDC00}.
\end{problem}

\section{Conclusions} \label{sec-conclu}
In this paper, the numerical computations of initial value system of two dimensional \emph{convection-diffusion} equations with both kinds of boundary conditions has been done by adopting three methods: modified exponential cubic B-splines DQM, modified trigonometric cubic B-splines DQM, and mECDQ method \cite{SPK16}, which transforms  the \emph{convection-diffusion} equation into a system of first order ordinary differential equations (ODEs), in time, which is solved by using SSP-RK54 scheme.

The methods are found stable for two space convection-diffusion equation by employing matrix stability analysis method. Section \ref{sec-num} shows that the proposed solutions are more accurate in comparison to the solutions by various existing schemes, and are in good agreement with the exact solutions.

 The order of accuracy of the proposed methods for the convection-diffusion problem with Dirichlet's boundary conditions is cubic whenever $\beta_x, \beta_y>0$ and otherwise it is super linear, in space. On the other hand, the order of accuracy of the proposed methods for the convection-diffusion problem with Neumann boundary condition is quadratic with respect to $L_2$ error norms, see Table \ref{tab2.4}.

%\section*{\textbf{Acknowledgement}} \noindent The authors would like to thanks the anonymous referees for their time, effort, and extensive comments on the revision of the paper which improved the quality and presentation of the paper.

%\begin{figure}
%\centering
%\includegraphics[height=6.0cm, width=10.95cm]{SSPRK54F}
%\caption{Stability region for SSP-RK54 scheme} \label{Stab-R}
%\end{figure}

\end{document}